# A COMPUTATIONAL APPROACH TO CLASSIFYING LOW RANK MODULAR CATEGORIES

DANIEL CREAMER

ABSTRACT. This paper introduces a computational approach to classifying low rank modular categories up to their modular data. The modular data of a modular category is a pair of matrices, $(S, T)$. Virtually all the numerical information of the category is contained within or derived from the modular data. The modular data satisfy a variety of criteria that Bruillard, Ng, Rowell, and Wang call the admissibility criteria. Of note is the Galois group of the $S$ matrix is an abelian group that acts faithfully on the columns of the eigenvalue matrix, $s = (\frac{S_{ij}}{S_{0j}})$. This gives an injection from $\text{Gal}(\mathbb{Q}(S), \mathbb{Q}) \to \text{Sym}_r$, where $r$ is the rank of the category. Our approach begins by listing all the possible abelian subgroups of $\text{Sym}_6$ and building all the possible modular data for each group. We run each set of modular data through a series of Gröbner basis calculations until we either find a contradiction or solve for the modular data.

The effectiveness of this approach is shown by the two main results. The first is a complete classification of rank 6 non-self-dual and non-integral modular tensor categories, specifically any rank 6 non-integral non-self-dual modular category is isomorphic to a tensor product of the Fibonacci category and $\mathbb{Z}_3$ or their Galois conjugates. The second is a partial classification of the subgroups of $\text{Sym}_6$ that give rise to self-dual non-integral modular tensor categories. Specifically, we show that the following groups have no associated modular category, $\langle (01234) \rangle$, $\langle (0123) \rangle$, $\langle (01)(23), (02)(13) \rangle$, $\langle (0123)(45) \rangle$, $\langle (012), (345) \rangle$, $\langle (01), (2345) \rangle$, $\langle (01)(2345) \rangle$, $\langle (01), (23)(45), (24)(35) \rangle$, $\langle (01)(23)(45) \rangle$, $\langle (01)(23)(45), (24)(35) \rangle$, or $\langle (01)(23), (23)(45) \rangle$. It is known that following groups do have categories associated with them, $\langle (012) \rangle$, $\langle (01)(23) \rangle$, $\langle (012)(345) \rangle$, $\langle (01)(23)(45), (02)(13) \rangle$, and $\langle (012345) \rangle$. It is unknown but conjectured that the following groups do not have a modular category associated to them, $\langle (01) \rangle$, $\langle (01), (234) \rangle$, and $\langle (01), (23)(45) \rangle$.

## 1. INTRODUCTION

The classification of fusion categories, and in particular modular tensor categories (MTC) has been an active research area for more than 20 years. Classifying fusion categories is a very broad problem. We have chosen here to focus on modular tensor categories, sometimes just called modular categories. To break this into a smaller problem in this paper we consider classifying MTC by their rank, $r$.

In 2016, Bruillard, Ng, Rowell, and Wang proved that there are finitely many MTC for a given rank[4]. In the same year they completed the classification of all modular categories through rank 5[3]. We have extended the classification with a partial classification of rank 6 and have created a computer assisted classification method to continue this classification program for low rank MTC.

A modular tensor category, $\mathcal{C}$, is a non-degenerate ribbon fusion category over $\mathbb{C}$. Let $\Pi_{\mathcal{C}}$ be the set of isomorphism classes of simple objects of the modular category $\mathcal{C}$. The rank of $\mathcal{C}$ is the finite number $r = |\Pi_{\mathcal{C}}|$. Some standard references for modular categories are Lectures on Tensor Categories [1], Tensor Categories [8], and Quantum Invariants of Knots and 3-Manifolds [12].

Every MTC has a pair of associated matrices, $(S, T)$, called the modular data. The modular data determine virtually all the numerical invariants of the category. The modular data are both $r \times r$ square matrices. The $S$ matrix can be used to immediately find both the dimension of the simple objects but also the Frobenius-Perron dimension of the simple objects (FPdim). Those dimensions can then be used to find the global dimension of the category and the Frobenius-Perron dimension of the category. Going farther, the $S$-matrix determines the fusion rules of the category. This is important because previous classifications of various fusion categories were done up to the fusion rules. Classifying all possible modular data is then at least as strong as classifying the possible fusion rules. Bruillard, Ng, Rowell, and Wang also determined a set of criteria that these two matrices must satisfy. They call it admissibility criteria. They conjecture that if two matrices are admissible then there is an associated MTC. [4][3]

The modular data also have a rich Galois structure. In particular,
$$\text{Gal}(\mathbb{Q}(S)/\mathbb{Q}) \leq \text{Gal}(\mathbb{Q}(T)/\mathbb{Q}).$$
The $T$ matrix is a diagonal matrix with roots of unity on the diagonal. Therefore, $G = \text{Gal}(\mathbb{Q}(S)/\mathbb{Q})$ is isomorphic to a subgroup of $(\mathbb{Z}/N\mathbb{Z})^\times$, where $N$ is the smallest positive integer such that, $T^N = \text{Id}$. Note that this group is therefore abelian. There is a related matrix, $s = (s_{ij}) = \left(\frac{S_{ij}}{S_{0j}}\right)$, which we call the eigenvalue matrix. It turns out that $G$ permutes the columns of $s$ faithfully. This gives an injection into the symmetric group on $r$ letters, $\text{Sym}_r$. Abelian subgroups of $\text{Sym}_r$ are well understood and easy to list. Furthermore, $G$ acts on the entries of $S$. This combined with some other symmetries of the $S$ matrix means that instead of 42 unknown entries of $S$ and $T$, there might actually only be 18.

In short, our method first determines a list of abelian subgroups of $\text{Sym}_r$. Then, for a given subgroup, we build the possible $S$ and $T$ matrices. We enter certain polynomial relations of the modular data into a Gröbner basis algorithm (GBA) and factor the result. We look for contradictions to the admissibility criteria and if we don't find any, refine the list of polynomials that are fed into the GBA in a variety of ways. We continue this until we have solved for the modular data or have proven the modular data is inadmissible.

We have two results. The first is a complete classification of non-integral, non-self-dual, modular tensor categories (NI-NSD-MTC) of rank 6. Namely, all NI-NSD-MTC of rank 6 are isomorphic to a tensor product of the Fibonacci MTC and $\mathbb{Z}_3$ MTC (or their Galois conjugates). The second is a partial classification of the abelian subgroups of $\text{Sym}_6$ that do have a non-integral, self-dual, modular tensor category (NI-SD-MTC) associated to them. Up to a relabeling of the simple objects the groups that are known to realize a MTC are $\langle(012)\rangle$, $\langle(01)(23)\rangle$, $\langle(012345)\rangle$, $\langle(012)(345)\rangle$, and $\langle(01)(23)(45), (02)(13)(45)\rangle$.

2. Previous Results

Modular categories have a lot of rich structure. But this makes them somewhat rare. Most of our classification method is eliminating possible modular data. We collect information that will do one of two things, either determine a polynomial relation to feed into the Gröbner basis algorithm or a fact about one the unknown constants that can be used to refine the Gröbner basis that is outputted.

**Definition 1.** *For a pair of matrices $(S,T)$ for which there exists a modular category with modular $(S,T)$, we will say $(S,T)$ is **realizable modular data**.[4]*

**Definition 2.** *Let $S, T \in GL_r(\mathbb{C})$ and define constants, $d_i := S_{0i}$, $\theta_i := T_{ii}$, $D^2 := \sum_i d_i^2$, and $p_\pm := \sum_i S_{0i}^2 \theta_i^{\pm 1}$. The pair $(S,T)$ is an **admissible modular data** of rank $\mathbf{r}$ if they satisfy the following conditions:*

(i) $d_i \in \mathbb{R}$, $S = S^t$, and $S\bar{S} = D^2 \text{Id}$. $T_{i,j} = \delta_{i,j}\theta_i$ with $N := \text{ord}(T) < \infty$
(ii) $(ST)^3 = p_+ S^2$, $p_+ p_- = D^2$, and $\frac{p_+}{p_-}$ is a root of unity
(iii) $N_{i,j}^k := \frac{1}{D^2} \sum_{a=0}^{r-1} \frac{S_{ia}S_{ja}\bar{S}_{ka}}{S_{0a}} \in \mathbb{N}$ for all $0 \leq i, j, k \leq (r-1)$
(iv) $\theta_i \theta_j S_{ij} = \sum_{k=0}^{r-1} N_{i^*j}^k d_k \theta_k$, where $i^*$ is the unique label such that $N_{ii^*}^0 = 1$.

(v) $\nu_n(k) := \frac{1}{D^2} \sum_{i,j=0}^{r-1} N_{ij}^k d_i d_j \left(\frac{\theta_i}{\theta_j}\right)^n$. Then $\nu_2(k) = 0$ if $k \neq k^*$ and $\nu_2(k) = \pm 1$ if $k = k^*$. Moreover, $\nu_n(k) \in \mathbb{Z}[e^{2\pi i/N}]$

(vi) $\mathbb{Q}(S) \subset \mathbb{Q}(T) = \mathbb{Q}(e^{2\pi i/N}) = \mathbb{Q}_N$, $\mathrm{Gal}(\mathbb{Q}(T)/\mathbb{Q})$ is isomorphic to an abelian subgroup of $\mathrm{Sym}_r$ and $\mathrm{Gal}(\mathbb{Q}(T)/\mathbb{Q}(S)) \cong (\mathbb{Z}/2\mathbb{Z})^l$ for some integer $l$.

(vii) The prime divisors of $D^2$ and $N$ coincide in $\mathbb{Z}[e^{2\pi i/N}]$

[4]

There are two facts that it is assumed the reader knows when looking at the list above. First the 0 label is unique and we label the $r$ isomorphism classes of simple objects from 0 to $r - 1$. The 0 object is often referred to as the vacuum because of the application of fusion categories to theoretical physics. For this paper, it is just important to remember that this is the only isomorphism class that cannot be relabeled. The second is the $*$ is an involution. It's map that takes a label to its dual label. For the purposes of the work in this paper, all we need to know is that it is an involution, namely $(i^*)^* = i$. Note that $N_{00}^0 = \frac{1}{D^2} \sum_a d_a^2 = 1 \implies 0^* = 0$.

A Gröbner basis algorithm takes a collection of polynomials and outputs a basis for the ideal of the polynomials that share the same solution as the initial collection over a given field. We use the field of rational numbers. Our approach is to factor this output and look at the parts of the basis that have distinct factors. For example, We might see $d_2(p^2 + D^2)$ as a polynomial in the outputted basis. From the admissibility criteria, $d_2 \neq 0$ (look at (iii)). Therefore, $p^2 + D^2$ must be 0. Now, we can add, $p^2 + D^2$ as a relation and rerun the GBA. Eventually we hope to either solve for the entries of $(S, T)$ or to find a contradiction.

The admissibility criteria give a foundation for building candidate pairs $(S, T)$, particularly, (i) and (ii). The first important fact, $S$ is self transpose. This immediately reduces the number of variables to solve for in the modular data. Not immediately clear in the admissibility criteria above, $S$ also has these two properties, $S_{ij} = S_{i^*j^*}$ and $\bar{S}_{ij} = S_{i^*j}$. Second, we get two polynomial relations, $S\bar{S} = D^2 \mathrm{Id}$ (**orthogonality relations**) and $pS^2 = (ST)^3$ (**twist relations**). These relations will be the foundations for the attempts to solve the modular data. In this paper, we do not use (iii) and (iv). These can be useful relations, but they introduce the variables, $N_{ij}^k$ that must be solved for. There is an important use for (v). It pairs with a theorem from Ng to give a powerful new relation. The theorem is at the moment unpublished but has been presented in at least two conferences.

**Theorem 1.** *If $p$ does not divide $N = \mathrm{ord}(T)$ for a modular category, $\mathcal{C}$, then there exists a simple object such that,*

$$X^2 = \sum_i \nu_p(X_i) X_i.$$

*Moreover, the label of $X$ in $\Pi_\mathcal{C}$ is in the orbit of 0. [10]*

Point (vi) is a product of the strong Galois symmetries at play. The Galois symmetries are possibly what makes this classification attempt feasible. The isomorphism from $\mathrm{Gal}(\mathbb{Q}(S)/\mathbb{Q})$ to some abelian subgroup of $\mathrm{Sym}_r$ is our starting point towards building modular data candidates. The restriction of $\mathrm{Gal}(\mathbb{Q}\mathbb{W}(T)/\mathbb{Q}(S))$ will give important information about $N$, the order of $T$. The prime divisors of $N$ are important later when considering the representation theory of a possible modular category.

**Theorem 2.** *Let $(S, T)$ be a realizable modular data. Then,*
*(a) $(S, T)$ is admissible.*
*(b) $(\sigma(S), \sigma(T))$ is realizable [4]*

This theorem is very useful, particularly (b). Part (a) says exactly what we expect, that the admissibility criteria is necessary (it may not be sufficient) for a pair of matrices to be modular data for some MTC. But part (b) allows for assumptions when classifying the data. It might be that $(\sigma(S), \sigma(T))$ has a desirable property that $(S, T)$ do not. A common property that appears in this manner is that for some $\sigma$, $\sigma(S_{0i}) \geq 1$ for all $i$. This says that for $S' = \sigma(S)$, $d'_i = S'_{0i} \geq 1$.

**Theorem 3.** *Let $G = \text{Gal}(\mathbb{Q}(S)/\mathbb{Q})$ in $\text{Sym}_r$ for some realizable $S$. Then,*
- *For every $\sigma \in G$, there is a sign function $\epsilon_\sigma : \Pi_\mathcal{C} \to \pm 1$ such that,*

$$S_{ij} = \epsilon_\sigma(\sigma(i))\epsilon_\sigma(j)S_{\sigma(i)\sigma^{-1}(j)}.$$

- *If $r$ is even, then $\prod_{i=0}^{r-1} \epsilon_\sigma(i) = (-1)^\sigma$.*

[11]

This is a slight abuse of notation. Throughout this paper we will use $G$ and its image in $\text{Sym}_r$ interchangeably. This makes reading lines like the first bullet point easier by not adding extra symbols. This theorem is the last fundamental result that makes this process feasible. Much like the fact that $S$ is self dual, the first bullet point reduces the number of variables of the $S$ matrix. It also potentially adds variables in the choices of the sign function. However, it turns out that the already known symmetries of the $S$-matrix combined with this bullet point put significant restrictions on the values of $\epsilon_\sigma(i)$ related to each other. A priori it appears that for rank 6, there are $2^6$ choices of the sign function. In practice there are usually less than 8 relevant choices of the sign function. This makes it more expedient to treat each choice of the sign function as its own modular data candidate. Doing so eliminates the need to use variables for $\epsilon_\sigma(i)$ and frequently sign choices lead to virtually identical contradictions.

**Definition 3.** *The **fusion matrices**, $N_i$, are defined to be $(N_i)_{kj} = N_{ij}^k$.*

**Definition 4.** *Let $s_{ij} = \frac{S_{ij}}{S_{0j}}$, then we call $s$ the **eigenvalue matrix**.*

The eigenvalue matrix is so named because it contains the eigenvalues of $N_i$ in the $i^{th}$ row, in particular this means that the entries of $s$ and $S$ are algebraic integers. This isn't immediately obvious until noticing that point $(iii)$ of definition 2 (called the Verlinde Formula) can be rewritten as $SN_iS^{-1} = D_i$ where $D_i = (\delta_{ab}\frac{S_{ia}}{S_{0a}})_{ab}$. This shows that the $S$ matrix diagonalizes $N_i$, and the diagonal is the $i^{th}$ row of the eigenvalue matrix. The eigenvalue matrix is also where the isomorphism into $\text{Sym}_r$. The isomorphism comes from the fact that $S$ simultaneously diagonalizes all of the fusion matrices and thus $\text{Gal}(\mathbb{Q}(S)/\mathbb{Q})$ permutes the columns of $s$ [2]. Another immediate fact from this formula is that $\mathbb{Q}(S)$ is a Galois extension of $\mathbb{Q}$. It should be clear that $\mathbb{Q}(S)$ and $\mathbb{Q}(s)$ are the same field. The entries of $s$ are fractions of the entries in $S$, and all the entries is $S$ are products of entries in $s$. But each row of $s$ is the entire set of eigenvalues of an integer matrix. Therefore, $\mathbb{Q}(S)$ is the splitting field of the characteristic polynomial of the integer matrices, $N_i$. Those polynomials are in $\mathbb{Z}[x]$, so the field extension is Galois.

We know that the fusion matrices have non-negative integer entries, i.e. they have a Frobenius-Perron eigenvalue, a positive real eigenvalue that is larger in modulus than any eigenvalue (real or complex). For fusion categories this is called the **FPdim** of $X_i$, where $X_i$ is the isomorphism class of simple object of label $i$. For a fusion category $\mathcal{C}$, the $\text{FPdim}(\mathcal{C} = \sum \text{FPdim}(X_i)^2$. One of the columns of $s$ will contain the FPdim of $X_i$ in the $i^{th}$ row [4]. We will call the corresponding column of $S$ the FPcol.

Since the entries in the FPcol are just a product of $d_i$ and the FPdims, then the FPcol is entirely negative or entirely positive depending on the sign of $d_i$. Then the orthogonality relations imply that there cannot be two columns of $S$ where each entry has the same sign. Also, since neither $d_i$ nor any FPdim can be 0, the FPcol has no 0's as entries.

I add a corollary to theorem 2.

**Corollary 4.** *The eigenvalue matrix of $\sigma(S)$ is $\sigma(s)$.*

*Proof:* Let $S' = \sigma(S)$ and $d'_i = S'_{0i}$. This follows immediately once we write $S_{ij}$ as $d_i s_{ij}$. Then, $\sigma(S_{ij}) = \sigma(d_i)\sigma(s_{ij})$. But, $d_i = S_{0i} \implies \sigma(d_i) = d'_i$. Thus it's clear that the eigenvalue matrix of $S'$ is $\sigma(s)$. $\square$

This corollary shows that $(\sigma(S), \sigma(T))$ have the same associated subgroup of $\text{Sym}_r$ as $(S,T)$. If the FPcol of $S$ is column $i$. Then for purposes of classification, we can instead consider $(\sigma(S), \sigma(T))$ and therefore assume the FPcol is any column with label in the orbit of $i$. In, $S$ the FPcol can often be shown to be the in the orbit 0. Therefore without loss of generality we assume the FPcol is column 0 and that $d_i \geq 1$ for all $i$.

**Definition:** A fusion category, $\mathcal{C}$ is to be

(i) **weakly integral** if $\text{FPdim}(\mathcal{C}) \in \mathbb{Z}$.
(ii) **integral** if $\text{FPdim}(X_i) \in \mathbb{Z}$ for all $i \in \Pi_\mathcal{C}$
(iii) **pointed** if $\text{FPdim}(X_i) = 1$ for all $i \in \Pi_\mathcal{C}$

**Theorem 5.** *Let $\mathcal{C}$ be a modular category. Then, $\mathcal{C}$ is integral if and only if $d_i \in \mathbb{Z}$ for all $i \in \Pi_\mathcal{C}$ if and only if $\mathcal{O}_0 = \{0\}$. [3]*

Weakly integral modular categories have been classified up through rank 7. In particular for rank 6 this means that the only abelian subgroups of $\text{Sym}_6$ to consider are the ones that move 0.

**Definition 5.** *A modular category is self dual if $i^* = i$ for all $i \in \Pi_\mathcal{C}$.*

Now look at the admissibility criteria again to see that this is equivalent to saying that $\mathcal{C}$ is self dual if and only if $S$ is a real matrix. Recall that the admissibility criteria imply that $\bar{S}_{ij} = S_{i^*j}$. So, if $\mathcal{C}$ is self dual then for all $i,j \in \Pi_\mathcal{C}$, $\bar{S}_{ij} = S_{ij} \implies S$ is real. This leads to an important fact about $G$, complex conjugation (which is always in the Galois group), corresponds switching dual labels. Therefore, in a non-self-dual category, there must exist an element of $G$ that switches dual labels.

**Theorem 6.** *For $i \in \Pi_\mathcal{C}$, if $j \in \mathcal{O}_i$ then $i$ and $j$ are both self-dual labels or are both non-self-dual labels. [3]*

This theorem gives a bit more structure to potential subgroups of $\text{Sym}_6$ when considering the non-self-dual cases. Recall that the admissibility criteria forces $0^* = 0$ and that in the non-integral cases, the 0 label is not fixed. This means for rank 6 there are 0, 1, or 2 pair of non-self-dual classes of simple objects.

**Theorem 7.** *If $\sigma^2 = \text{Id}$ and $\epsilon_\sigma$ is the corresponding sign function, $\sigma$ fixes at least one element of $\Pi_\mathcal{C}$, then $\epsilon_\sigma(i) = \epsilon_\sigma(\sigma(i))$.[3]*

We've chosen to use notation from [11], $S_{ij} = \epsilon_\sigma(\sigma(i))\epsilon_\sigma(j)S_{\sigma(i)\sigma^{-1}(j)}$. But, [3] uses $S_{ij} = \epsilon_\sigma(i)\epsilon_{\sigma^{-1}}(j)S_{\sigma(i)\sigma^{-1}(j)}$. These $\epsilon_\sigma$ are defined slightly differently. But this conclusion still holds. To see this, let $\alpha$ be the sign function that corresponds to $\sigma$ under [11] and $\beta/\beta^{-1}$ be the sign functions that correspond to $\sigma$ and $\sigma^{-1}$ respectively from [3]. Now, since $\sigma$ is assumed to be an involution, $\beta = \beta^{-1}$. Let $j$ be fixed by $\sigma$ and $i \in \Pi_\mathcal{C}$. Observe,

$$S_{0j} = \alpha(\sigma(0))\alpha(j)S_{\sigma(0)j} = \alpha(\sigma(0))\alpha(j)\alpha(\sigma^2(0))\alpha(j)S_{0j}$$
$$\implies \alpha(0)\alpha(\sigma(0)) = 1 \implies \alpha(0) = \alpha(\sigma(0))$$
$$S_{0i} = \alpha(\sigma(0))\alpha(i)S_{\sigma(0)\sigma^{-1}(i)} = \beta(0)\beta(i)S_{\sigma(0)\sigma^{-1}(i)}$$
$$\implies \alpha(0)\alpha(i) = \beta(0)\beta(i)$$
$$S_{i0} = \alpha(\sigma(i))\alpha(0)S_{\sigma(i)\sigma^{-1}(0)} = \beta(i)\beta(0)S_{\sigma(i)\sigma^{-1}(0)}$$
$$\implies \alpha(\sigma(i))\alpha(0) = \beta(i)\beta(0) \implies \alpha(i) = \alpha(\sigma(i))$$

Recall that a priori there are $2^6$ possibilities for the sign function corresponding to $\sigma$. We should first note, that for the symmetries of the $S$ matrix, we care about $\epsilon_\sigma(i)\epsilon_\sigma(j)$ for all $i, j \in \Pi_\mathcal{C}$. This means that for computations we can assume $\epsilon_0 = 1$. But that still leaves $2^5$ possibilities. Theorems like this are helpful in reducing that even farther. Then, much like the first line in the above calculation we can get similar relations to continue to reduce the number of possibilities.

## 3. Preliminary Results

This section contains original work that again primarily serves the same two functions as the facts from the previous section.

**Definition:** The *s*-**polynomial** of $s_{ij}$ is $\prod_{k \in \mathcal{O}(j)}(x - s_{ik})$, where $\mathcal{O}(j)$ is the orbit of $j$ under the permutation action of $G$ on the columns of $s$.

Note that $d_i \in \text{Gal}(S/\mathbb{Q})$ and $d_i \in \text{Gal}(s/\mathbb{Q}) \implies \text{Gal}(S/\mathbb{Q}) = \text{Gal}(s/\mathbb{Q})$.

**Theorem 8.** *The s-polynomial has integer coefficients and is a power of the minimal polynomial of $s_{ij}$.*

*Proof:* This is clear. The Galois group clearly permutes the linear factors of the *s*-polynomials. Therefore the polynomial itself is fixed by the entire group, i.e. the coefficients are rational numbers. They are also sums of products of algebraic integers and therefore algebraic integers themselves. Thus, the coefficients are integers. Note that the leading coefficient is necessarily a one. Therefore all irreducible factors over $\mathbb{Z}[x]$ also have a leading coefficient of one, i.e. they are all minimal polynomials of some algebraic integer in the given orbit. But each term in the given orbit necessarily have the same minimal polynomial. The only possible irreducible factor is the minimal polynomial of $\frac{S_{ij}}{S_{0j}}$. □

These *s*-polynomials introduce new variables and relations for the Gröbner basis calculations. Specifically, the variables are the coefficients of the *s*-polynomials. The new variables mean a more complicated polynomial ideal, but it also introduces new information to be used when looking at the factored basis. Namely, that the new variables are integers.

**Theorem 9.** *If $\sigma = (01)$ is in the Galois group, and the rank, $r$, of the MTC is at least five, then*
*i)* $\frac{1}{d_1} + d_1$, $D^2/d_1$ *and* $d_i^2/d_1$ *are rational integers for $i \geq 2$.*
*ii) There exist $i, j \geq 2$ such that $\epsilon_i = -\epsilon_j$, and in this case $S_{ij} = 0$.*
*iii) $d_1 > 0$*
*iv) $\theta_1 = 1$*

*Proof:* This is a generalization of Lemma 3.6 from [3]. The proof is nearly identical, but by not requiring the Galois group to be generated by $(01)$, we no longer immediately force the FPcol to be $S_0$ or $S_1$ where $S_k$ is the $k^{th}$ column of $S$, although it does still turn out to be the case.

First we must show that the orbit of 0 is in fact $\{0, 1\}$. Suppose $\tau(0) \neq 0$ or 1. Then, $\tau\sigma\tau^{-1} = (\tau(0)\tau(1)) \neq (01)$. But then $\tau \notin G$ since $G$ is abelian.

By 3.5 of [3] we see that $S_{11} = 1$. Then by the s-polynomial of $d_1$, $(d_1 + \frac{1}{d_1}) \in \mathbb{Z}$. Consider the s-polynomial of $d_j$ for $j \geq 2$. This shows that $\frac{d_j^2}{d_1} \in \mathbb{Z}$. Then, $\frac{D^2}{d_1} = (\frac{1}{d_1} + d_1) + \sum_{j \geq 2} \frac{d_j^2}{d_1} \in \mathbb{Z}$.

For ii) note that $S_{1j} = \epsilon_0 \epsilon_j d_j$ and apply orthogonality of the columns to get, $2d_1 = \sum_{j=2}^{r-1} \epsilon_0 \epsilon_j d_j^2$, divide by $d_1$ and get $2 = \sum_{j=2}^{r-1} \epsilon_0 \epsilon_j d_j^2 / d_1$. Note that $d_j^2/d_1$ has the same sign for all $j \geq 2$. Since the rank is at least 5 and we can't have 2 be the sum of 3 or more integers of the same sign, we see there must be some $i, j \geq 2$ such that, $\epsilon_i = -\epsilon_j$.

But now, observe that $S_{ij} = \epsilon_i \epsilon_j S_{ij} = -S_{ij} \implies S_{ij} = 0$.

iii) An immediate consequence is that every column except $S_0$ and $S_1$ must have a 0 entry. Therefore either $S_0$ or $S_1$ is the FPcol and is strictly positive, i.e. $d_1 > 0$ since it's in both columns (a priori $S_1$ could be strictly negative but since 1 is in $S_1$ we know it can't be strictly negative).

iv) Take, the twist equation (for simplicity let $p = p^+$), $0 = pS_{ij}$ where $i, j$ are the indexes given in ii). Then, we see that,

$$0 = pS_{ij} = \theta_i \theta_j (d_i d_j - \theta_1 d_i d_j + 0) = \theta_i \theta_j d_i d_j (1 - \theta_1)$$

This implies that $\theta_1 = 1$. $\square$

This theorem allows many groups to be tackled simultaneously. By only looking at the structure of $S$ given by (01) and the choices on the corresponding sign function $\epsilon$, we can potentially eliminate all groups that contain (01) at once. This would be difficult without knowing more about the $S$ matrix, except the choices on the sign function provide 0's as entries in the $S$ matrix and greatly simplify many of the polynomial relations found in the orthogonality and twist relations.

**Corollary 10.** *Let* $(01) \in G$. *Let* $\sigma = (01)$.
  *i) If* $\epsilon_i \neq \epsilon_j$ *for all* $j \neq i$, $2 \leq j \leq r-1$, *then* $\theta_i^2 + \theta_i + 1 = 0$.
  *ii) If* $\epsilon_i \neq \epsilon_j$ *for all* $j \neq i, i^*$, $2 \leq j \leq r-1$ *and* $i \neq i^*$, *then* $\theta_i^2 + 1 = 0$.

*Proof:* i) Let's assume $\epsilon_2 \neq \epsilon_j$ for $j > 2$. By theorem 3, $S_{2j} = 0$. Let $zd_2 = S_{22}$. We know that $z$ is an integer because $z = \frac{S_{22}}{d_2}$ is fixed by all elements of the Galois group (or is possibly moved to the 0's that populate the rest of the row and is therefore 0, itself an integer). Let $S_k$ be the $k^{th}$ column of $S$. Then $S_0 \cdot S_2 = 0 = d_2(1 + \epsilon_1 \epsilon_2 d_1 + zd_2) \implies 1 + \epsilon_2 d_1 + zd_2 = 0$.

$pS_{ij} = \theta_i \theta_j \sum_{k=0}^{r-1} S_{ki} S_{kj} \theta_k$. Then, using $pS_{00}, pS_{01}, pS_{02}, pS_{22}$ and the relation, $1 + \epsilon_1 \epsilon_2 d_1 + zd_2$, we can deduce,

$$pz^2 = \theta_2 z^2 (\theta_2 - 1)$$
$$pz^2 = z^2 + 2\theta_2$$

Which gives,

$$0 = z^2 \theta_2^2 + \theta_2(2 - z^2) + z^2(1)$$

an integer polynomial of $\theta_2$, i.e. the minimal polynomial, $f$, of $\theta_2$ divides it. A quick check shows that $\theta_2$ isn't $\pm 1$. If $\theta_2 = 1 \implies p = 0$. If $\theta_2 = -1 \implies p = 2 \implies D^2 = 4 < 5$. But we can assume $D^2 \geq 5$ because we can assume the FPcol is column 0 and therefore $d_i \geq 1$ for all $i$.

$$f = x^2 + x \frac{2 - z^2}{z^2} + 1$$

is the minimal polynomial of $\theta_2 \implies \frac{2-z^2}{z^2} = 0, \pm 1$. Both, 0 and -1 give contradictions. If 0, then $z$ isn't an integer. If -1, then $2 = 0$. So, $\theta_2$ is a third root of unity.

ii) Assume for simplicity, $2^* = 3$ and $\epsilon_2 \neq \epsilon_i$ for all $i > 3$. Let $\epsilon = \epsilon_0 \epsilon_2$. Note that $0 = S_0 \cdot S_2 = d_2(1 + \epsilon d_1 + a_1 + a_2) \implies 1 + \epsilon d_1 = -(a_1 + a_2)$. Then,

$$pS_{20} = pd_2 = \theta_2(d_2 + \epsilon d_1 d_2 + \theta_2 d_2 a_1 + \theta_2 d_2 a_2)$$
$$\implies p = \theta_2(a_1 + a_2)(\theta_2 - 1).$$
$$\implies p(a_1 + a_2) = \theta_2(a_1 + a_2)^2(\theta_2 - 1)$$
$$pS_{22} = pa_1 = \theta_2^2(d_2^2 + d_2^2 + \theta_2 a_1^2 + \theta_2 a_2^2)$$
$$pS_{23} = pa_2 = \theta_2^2(d_2^2 + d_2^2 + \theta_2 a_1 a_2 + \theta_2 a_2 a_1)$$
$$\implies p(a_1 + a_2) = \theta_2^2(4d_2^2 + \theta_2(a_1 + a_2)^2)$$

But then,

$$p(a_1 + a_2) = p(a_1 + a_2) \implies$$
$$\theta_2(a_1 + a_2)^2(\theta_2 - 1) = \theta_2^2(4d_2^2 + \theta_2(a_1 + a_2)^2)$$
$$(a_1 + a_2)^2(\theta_2 - 1) = \theta_2(4d_2^2 + \theta_2(a_1 + a_2)^2)$$

Let $z = \frac{(a_1 + a_2)^2}{d_2^2}$. Arrange as a polynomial of $\theta_2$ and divide through by $d_2^2$:

$$z\theta_2^2 + \theta_2(4 - z) + z = 0$$

Observe that $\frac{a_1 + a_2}{d_2}$ is an integer by the $s$-polynomial of $\frac{a_1}{d_2}$ (recall that the rest of row 2 after column 3 is entirely filled with 0's and therefore cannot be in the orbit of column 2). Therefore $z$ is an integer and a perfect square. This also forces the polynomial above to have integer coefficients, i.e. the minimal polynomial of $\theta_2$, $\min(\theta_2)$, must divide $zx^2 + x(4 - z) + z$. We immediately know that the degree of $\min(\theta_2)$ is either 1 or 2. If the degree is 1, then $\theta_2 = \pm 1$. If $\theta_2 = 1$, then, $0 = p = \theta_2(a_1+a_2)(\theta_2-1)$. If $\theta_2 = -1$, then $z\theta_2^2+\theta_2(4-z)+z = 0 \implies z+z-4+z = 0 \implies 3z = 4$. But then $z$ is not an integer.

Since we can conclude that the degree of $\min(\theta_2) = 2$, then $zx^2 + x(4 - z) + z = z * \min(\theta_2)$, i.e. $\min(\theta_2) = x^2 + x\frac{4-z}{z} + 1$. Now, $\theta_2$ is a root of unity, and all minimal polynomials of degree 2 of roots of unity are known, we say that $\frac{4-z}{z} = 0, \pm 1$. If $\frac{4-z}{z} = 1 \implies z = 2$. But 2 isn't a perfect square. If $\frac{4-z}{z} = -1 \implies 4 = 0$. If $\frac{4-z}{z} = 0 \implies z = 4 \implies \min(\theta_2) = \theta_2^2 + 1$. $\square$

This corollary, provides some additional information but about the $T$ matrix. Such information about $T$ matrix very frequently leads to a near solution of the modular data or an immediate contradiction.

**Theorem 11.** *Let $A = \{a \mid X_a \text{ is SD}\}$. Let $B = \{b \mid X_b \text{ is NSD}\}$. Let $\epsilon$ be the sign function for $\sigma$ and $\delta$ be the sign function for $\tau$.*

*i) For all $\sigma \in G$, $\epsilon_b = \epsilon_{b^*}$ for all $b \in B$.*

*ii) Furthermore, if for all $a \in A$, $\tau(a) = \sigma(a)$ and for all $b \in B$, $\tau(b) = \sigma(b)$ or $\sigma(b)^*$ then $\epsilon_i \epsilon_j = \delta_i \delta_j$ for all $i, j$. If this is true we say $\epsilon = \delta$.*

*Proof:* i) First note that $\sigma(b^*) = \sigma(b)^*$. This is easy to observe since $G$ is abelian and complex conjugation is always in the Galois group. Let $\gamma$ be complex conjugation, i.e. $\gamma(i) = i^*$. Then, $\sigma(b^*) = \sigma\gamma(b) = \gamma\sigma(b) = \sigma(b)^*$. Combine this with, $S_{0b} = S_{0b^*} \ \forall \ b \in \{0, 1, ..., r-1\}$, and then observe:

$$\epsilon_{\sigma(0)} \epsilon_b S_{\sigma(0)\sigma^{-1}(b)} = S_{0b} = S_{0b^*} = \epsilon_{\sigma(0)} \epsilon_{b^*} S_{\sigma(0)\sigma^{-1}(b^*)}$$
$$= \epsilon_{\sigma(0)} \epsilon_{b^*} S_{\sigma(0)\sigma^{-1}(b)^*} = \epsilon_{\sigma(0)} \epsilon_{b^*} S_{\sigma(0)\sigma^{-1}(b)}$$

Since $S_{0b} \neq 0$, $S_{\sigma(0)\sigma^{-1}(b)} \neq 0 \implies \epsilon_b = \epsilon_{b^*}$.

ii) Now assume $\sigma$ and $\tau$ as in part ii). But then, $\sigma^{-1}$ and $\tau^{-1}$ must share the same assumptions. Note that it's enough to show $\epsilon_{\sigma(0)}\epsilon_i = \delta_{\sigma(0)}\delta_j$ for all $j$. Because what this really means is that all possible combination of products agree. But that's the same as finding out which $j$'s agree with a fixed $i$. In this case the fixed $i$ is $\sigma(0) = \tau(0)$.

So, observe
$$S_{0a} = \epsilon_{\sigma(0)}\epsilon_a S_{\sigma(0)\sigma^{-1}(a)} = \delta_{\sigma(0)}\delta_a S_{\tau(0)\tau^{-1}(a)}$$
$$= \delta_{\sigma(0)}\delta_a S_{\sigma(0)\sigma^{-1}(a)} \implies \epsilon_{\sigma(0)}\epsilon_a = \delta_{\tau(0)}\delta_a$$

because $S_{0a} \neq 0$.

This result holds clearly then if $\sigma(b) = \tau(b)$. So assume $\sigma(b) = \tau(b)^*$.

Similarly, $S_{0b} = \epsilon_{\sigma(0)}\epsilon_b S_{\sigma(0)\sigma^{-1}(b)} = \delta_{\sigma(0)}\delta_b S_{\tau(0)\tau^{-1}(b)} = \delta_{\sigma(0)}\delta_b S_{\sigma(0)\sigma^{-1}(b)^*} = \delta_{\sigma(0)}\delta_b S_{\sigma(0)\sigma^{-1}(b)}$
$\implies \epsilon_{\sigma(0)}\epsilon_a = \delta_{\tau(0)}\delta_b$ because $S_{0b} \neq 0$. $\square$

The first point further restricts the possible choices for the sign functions in the non-self-dual case. It continues a pattern for objects related to dual labels, $d_i = d_{i^*}$, $\theta_i = \theta_{i^*}$, and now $\epsilon(i) = \epsilon(i^*)$.

The second point's use is a bit more subtle. It also gives more restrictions on $\epsilon$ because now $\epsilon$ must satisfy $S_{ij} = \epsilon(\sigma(i))\epsilon(j)S_{\sigma(i)\sigma^{-1}(j)}$ for two slightly different $\sigma$. Under certain sign choices this may yet again force certain entries of the $S$ matrix to be 0.

**Theorem 12** (C, Munoz). *Let $\sigma = (012) \in G$. Let $r \geq 5$*

*i) There are four sign choices,*

*1) $\epsilon_0 = \epsilon_1 = \epsilon_2 = \epsilon_3$  2) $\epsilon_0 = \epsilon_1 \neq \epsilon_2 = \epsilon_3$  3) $\epsilon_0 = \epsilon_2 \neq \epsilon_1 = \epsilon_3$  4) $\epsilon_0 = \epsilon_3 \neq \epsilon_1 = \epsilon_2$*

*ii) When classifying modular data, we may assume sign choice 3.*

*Proof:* i) Observe, if $i \geq 3$, $S_{0i} = \epsilon_1\epsilon_i S_{1i} = \epsilon_1\epsilon_2 S_{2i} = \epsilon_0\epsilon_1\epsilon_2\epsilon_i S_{0i} \implies \epsilon_0\epsilon_1\epsilon_2\epsilon_i = 1$ and $\epsilon_3 = \epsilon_j$ for all $j \geq 3$. This leaves the 4 sign choices above.

ii) Assuming sign choice 1, there are no sign changes throughout the $S$-matrix as variables are moved around by $\sigma$. This means for every column there are at least the other columns with the exact same entries but possibly in different rows. Therefore, each set of three columns would share the same signs of the corresponding entries. Consider the FPcol, it's entirely negative or entirely positive, and so would be two other columns, this is already stated to be a contradiction.

Let $j > 2$, then as shown above $S_{0j} = d_j$, $S_{1j} = \epsilon_1\epsilon_3$, $S_{2j} = \epsilon_1\epsilon_2$. Recall $\epsilon_0\epsilon_1\epsilon_2\epsilon_3 = 1$. Now in sign choices 2 through 4 at least one of $S_{1j}$ and $S_{2j}$ will be negative. This means that for every column $j$ where $j > 2$, the column contains both $d_j$ and $-d_j$ and cannot be all negative or all positive. Thus no such column can be the FPcol. The FPcol must be either the 0, 1, or 2 column. Those labels share an orbit. For the purposes of classifying the modular data, we can assume the FPcol is column 0 and therefore that $d_i \geq 1$ for all $i$.

Under sign choice 2, $C_0 \cdot C_1 \implies d_1 + d_2 - d_1d_2 = \sum_{i=3}^{r-1} d_i^2 \geq 2$. But the surface, $d_1 + d_2 - d_1d_2$ has a maximum value of 1 under the restriction, $d_1, d_2 \geq 1$. Therefore both sign choices 1 and 2 are impossible.

Finally, sign choices 3 and 4 are actually relabelings of each other, via switching labels 1 and 2. Let $S3$ be the matrix of sign choice 3 and similarly $S4$ be the matrix of sign choice 4.

First if $i \geq 3$, $S3_{i1} = S4_{i2} = d_i$ and $S3_{i2} = S4_{i1} = -d_i$. Now we need to check that $S3_{12} = S4_{12}$, $S3_{22} = S4_{11}$, and $S3_{22} = S4_{11}$ remembering that in S4 we switch $d_1$ and $d_2$. Indeed, $S3_{12} = S4_{12} = -1$, $S3_{22} = d_1$ and $S4_{11} = d_2$, and $S3_{11} = -d_2$ and $S4_{22} = -d_1$. This shows that $S3$ and $S4$ are the same matrix up to switching labels 1 and 2. $\square$

This theorem may appear to have limited use, but it does apply to all ranks $\geq 5$.

Per the previous theorems, the Galois group of a MTC is an abelian subgroup of $\text{Sym}_6$, recall that the labeling starts at 0 and goes to 5. The orbit of self dual labels is entirely self dual and therefore the orbit of NSD labels are entirely NSD. The trivial object must always be labeled by 0 and is self dual. But, we are free to choose a labeling for the other objects. In a NSD rank 6 category there can be one pair or two pair of NSD objects.

**Theorem 13.** *1) If $\mathcal{C}$ is a rank 6 NI-NSD-MTC with two pair of simple objects, then it has one of the following Galois Groups (up to relabeling):* $\langle(01),(23),(45)\rangle$, $\langle(01)(23),(23)(45)\rangle$, $\langle(01),(23)(45)\rangle$, $\langle(01)(2435)\rangle$, $\langle(01),(2435)\rangle$, $\langle(01),(23)(45),(24)(35)\rangle$, *and* $\langle(01)(24)(35),(23)(45)\rangle$.

*2) If $\mathcal{C}$ is a NI-NSD-MTC with one pair of NSD objects, then it has one of the following Galois Groups (up to relabeling):* $\langle(01),(23),(45)\rangle$, $\langle(01)(23),(45)\rangle$, $\langle(012),(45)\rangle$, $\langle(01),(45)\rangle$, $\langle(01)(23),(02)(13),(45)\rangle$, *and* $\langle(0213),(45)\rangle$.

If the are two pair of NSD objects, let $2^* = 3$ and $4^* = 5$. If there is one pair of NSD objects let $4^* = 5$. We know that NSD categories necessarily have complex entries in the $S$ matrix in the columns with NSD labels. We also know from the $S$ matrix symmetries that complex conjugation switches columns. Therefore, $(23),(45) \in G$ or $(45) \in G$. Since the MTC is NI, the 0 label is not fixed. I will choose the orbit of 0 to be the first $k$ labels, where $k$ is the size of the orbit. Under these restrictions the possible Galois groups are the ones listed or are relabelings of the non-trivial SD objects. □

**Theorem 14.** *If $\mathcal{C}$ is a rank 6 NI-SD-MTC, then it has one of the following Galois Groups (up to relabeling):* $\langle(01234)\rangle$, $\langle(012345)\rangle$, $\langle(0123)\rangle$, $\langle(01)(23),(02)(13)\rangle$, $\langle(0123)(45)\rangle$, $\langle(01)(23)(45),(02)(13)\rangle$, $\langle(012)\rangle$,
$\langle(012),(345)\rangle$, $\langle(012)(345)\rangle$, $\langle(01),(2345)\rangle$, $\langle(01)(2345)\rangle$,
$\langle(01),(23)(45),(24)(35)\rangle$, $\langle(01)(23)(45),(24)(35)\rangle$, $\langle(01),(234)\rangle$, $\langle(01)(23)\rangle$, $\langle(01),(23)(45)\rangle$, $\langle(01)(23)(45)\rangle$, $\langle(01)(23),(23)(45)\rangle$, *or* $\langle(01)\rangle$.

There is only one group each where the size of the orbit of the 0 label is 5 or 6, $G = \langle(01234)\rangle$ or $\langle(012345)\rangle$. If the size of the orbit of the 0 label is 4, then restricted to the orbit of the 0 label, $G = \langle(0123)\rangle$ or $\langle(01)(23),(02)(13)\rangle$. From an earlier theorem, $(45) \notin G$. So, when viewing the action of $G$ on all of the labels, $G = \langle(0123)(45)\rangle$ or $\langle(01)(23)(45),(02)(13)\rangle$.

If the size of the orbit of the 0 label is 3, and the size of the orbit of the 3 label is 1 or 3, $G = \langle(012)\rangle$, $\langle(012),(345)\rangle$, or $\langle(012)(345)\rangle$. All groups with an orbit of size 3 for the 0 label and 2 for the 3 label will have a two-cycle that fixes 0. Therefore there are none.

Assume the orbit of the 0 label is $\{0,1\}$. If the size of the orbit of the 2 label is 4, then $G = \langle(01),(2345)\rangle$, $\langle(01)(2345)\rangle$, $\langle(01),(23)(45),(24)(35)\rangle$, or $\langle(01)(23)(45),(24)(35)\rangle$. If the size of the 2 label is 3, then $G = \langle(01),(234)\rangle$. If the the size of the obrit of the 2 label is 2, then $G = \langle(01),(23)(45)\rangle$, $\langle(01)(23)\rangle$, $\langle(01)(23)(45)\rangle$, or $\langle(01)(23),(23)(45)\rangle$. Finally that leaves $G = \langle(01)\rangle$.
□

A quick note about the relabeling, because we've fixed these groups, it's important that any relabeling we do preserves the group structure. If there is a relabeling that preserves the group structure, there will often be two sign choices that produce equivalent modular data under that relabeling.

## 4. Non-Integral, Non-Self-Dual, Modular Categories

The remaining sections are heavy on computation. Here are some notes to help understand the presented notes. In order to avoid confusion when looking at the factored Gröbner bases, I use $D$ for $D^2$. In all groups, $\epsilon$ is the corresponding sign function for the indicated $\sigma$. If necessary, $\delta$ is the corresponding sign function for the indicated $\tau$. Also, $\epsilon_i = \epsilon(i)$ to visually simplify the results. I will also $t_i = \theta_i$ to make the code easier to read and copy over. The first Gröbner basis run (unless otherwise stated) is run with the ideal generated by the orthogonality and the twist relations. The files for each computation and the output are available upon request. The code is nearly identical for each group. The only differences being the slight changes in the initial data (the ring, the S/T matrices, and occasionally an extra relation for the first Gröbner calculation)

Figure 1. Sample Code

```
|-- S1, 4^* = 5 and all other indices are self dual

R = QQ[D,d1,d2,d3,d4,a1,a2,a3,b1,b2,p,t2,t3,t4,e];

S = matrix {{1,d1,d2,d3,d4,d4},{d1,1,e*d2,e*d3,-e*d4,-e*d4},{d2,e*d2,a1,a2,0,0},
{d3,e*d3,a2,a3,0,0},{d4,-e*d4,0,0,a4,a5},{d4,-e*d4,0,0,a5,a4}};
A = matrix {{1,d1,d2,d3,d4,d4},{d1,1,e*d2,e*d3,-e*d4,-e*d4},{d2,e*d2,a1,a2,0,0},
{d3,e*d3,a2,a3,0,0},{d4,-e*d4,0,0,a5,a4},{d4,-e*d4,0,0,a4,a5}};

W = matrix {{D,0,0,0,0,0},{0,D,0,0,0,0},{0,0,D,0,0,0},{0,0,0,D,0,0},{0,0,0,0,D,0},
{0,0,0,0,0,D}};

T = matrix {{1,0,0,0,0,0},{0,1,0,0,0,0},{0,0,t2,0,0,0},{0,0,0,t3,0,0},
{0,0,0,0,t4,0},{0,0,0,0,0,t4}};

f1 = e^2-1;
f2 = t4^2+1;

I1 = ideal(S*A-W,p*S-T*S*T*S*T,f1,f2);

gbTrace = 3;
time G1 = gb(I1);

R = (flatten(entries(gens(G1))));

F = for r in R list factor(r);

toList(F_0)

for p in F do (
    pp = toList(p);
    for ppi in pp do (
        "S2-1" << "& " ;
        "S2-1" << toString(ppi) ;
        "S2-1" << " \n  \n &" ;
        );
    "S2-1" << " \n \n \n \n";
    )
"S2-1" << close

-- => pD = 0
```

and what new polynomial relations taken from the previous Gröbner basis run. All the code is run a Macbook Pro using Macaulay2 v 1.8.2. A sample can be seen above.

We first initialize the ploynomial ring over the rational numbers. Then we define $S$. The $A$ matrix is $\bar{S}$. Thus, $S * A = W$ is, $S\bar{S} = D^2 \text{Id}$, the orthogonality relations. We definte $T$. In this sample, we included some additional relations, $f1$ and $f2$, that are included in the first GBA. These were relations we were able to deduce before running the GBA the first time. Next comes the definition of the first ideal and the first Grobner calculation. After the Grobner calculation we have the program output the factored bases in a visually readable format. In particular we added a flag that can be searched for when a polynomial has at least two factors. For every factor after the first one, && is outputted. The polynomials that factor are the most common source of contradictions, as it happened in the sample code. But, they can give new relations too. For example we might see smething like $d_1(a_1 + a_2)$ appear in the output of our GBA. We know that $d_1$ is non-zero, therefore in the next run of the GBA we can add $a_1 + a_2$ as a relation and refine our ideal. In the case that we have new relations, we label each $h_k$ for some integer $k$, starting at 1. Then a second ideal would be constructed, $I2 = I1 + \text{ideal}(h_1 \ldots h_k)$, and the GBA run again. Sometimes the Gröbner basis calculation is halted at a particular degree. This gives a partial list of polynomials in the Gröbner ideal. This is done to find relations that will shorten the algorithm on successive runs. At the end of this section we've included tables that detail some of the steps used in running the code. The corresponding set of tables for section 5 is omitted for space purposes, but are available upon request.

**Theorem 15.** *There are no rank 6 NSD-NI-MTC's with $(01) \in G$.*

*Proof* Given theorems 9 and 10, we know that if $\epsilon_i \neq \epsilon_j$ for $2 \leq i, j \leq 5$ then $S_{ij} = 0$. We also know that at least one pair must exist. The three $S$-matrices below represent all the possible choices (up to relabelings) for $\epsilon_i \neq \epsilon_j$ given the MTC is NSD (recall $\epsilon_i = \epsilon_{i^*}$). In all three cases we also know that $t_1 = 1$.

$$S1 = \begin{pmatrix} 1 & d_1 & d_2 & d_2 & d_4 & d_4 \\ d_1 & 1 & d_2 & d_2 & -d_4 & -d_4 \\ d_2 & d_2 & a_1 & a_2 & 0 & 0 \\ d_2 & d_2 & a_2 & a_1 & 0 & 0 \\ d_4 & -d_4 & 0 & 0 & b_1 & b_2 \\ d_4 & -d_4 & 0 & 0 & b_2 & b_1 \end{pmatrix}$$

$S1$ is the only possible matrix for two pair of NSD objects. We can assume that both $t_2$ and $t_4$ are primitive 4th roots of unity, i.e. $t_2^2 + 1 = t_4^2 + 1 = 0$.

$$S2 = \begin{pmatrix} 1 & d_1 & d_2 & d_3 & d_4 & d_4 \\ d_1 & 1 & \epsilon_0\epsilon_2 d_2 & \epsilon_0\epsilon_3 d_3 & \epsilon_0\epsilon_4 d_4 & \epsilon_0\epsilon_4 d_4 \\ d_2 & \epsilon_0\epsilon_2 d_2 & a_1 & a_2 & 0 & 0 \\ d_3 & \epsilon_0\epsilon_3 d_3 & a_2 & a_3 & 0 & 0 \\ d_4 & \epsilon_0\epsilon_4 d_4 & 0 & 0 & b_1 & b_2 \\ d_4 & \epsilon_0\epsilon_4 d_4 & 0 & 0 & b_2 & b_1 \end{pmatrix}$$

$S_2$ assumes that $4^* = 5$ and that $\epsilon_2 = \epsilon_3 \neq \epsilon_4$. We can then assume that $t_4^2 + 1 = 0$.

$$S3 = \begin{pmatrix} 1 & d_1 & d_2 & d_3 & d_4 & d_4 \\ d_1 & 1 & \epsilon_0\epsilon_2 d_2 & \epsilon_0\epsilon_3 d_3 & \epsilon_0\epsilon_4 d_4 & \epsilon_0\epsilon_4 d_4 \\ d_2 & \epsilon_0\epsilon_2 d_2 & z_1 d_2 & 0 & 0 & 0 \\ d_3 & \epsilon_0\epsilon_3 d_3 & 0 & z_2 d_3 & z_3 d_4 & z_3 d_4 \\ d_4 & \epsilon_0\epsilon_4 d_4 & 0 & z_4 d_3 & b_1 & b_2 \\ d_4 & \epsilon_0\epsilon_4 d_4 & 0 & z_4 d_4 & b_2 & b_1 \end{pmatrix}$$

$S_3$ assumes that $4^* = 5$ and that $\epsilon_2 \neq \epsilon_3 = \epsilon_4$. We can then assume that $t_2^2 + t_2 + 1 = 0$. Observe, $S_{43} = S_{4^*3} = S_{53}$ and $S_{43} = S_{34}$. Therefore, $z_3 d_4 = z_4 d_3$.

Table 1 gives the details of each GBA calculation for $S-1$, $S-2$, and $S-3$. Let, $e = \epsilon_0\epsilon_2$.

Note that since columns 2-5 all have a 0, we can assume the FPcol is the 0 column and therefore, $d_i \geq 1$ for all $i \in \Pi_\mathcal{C}$.

Deduction 1: If $z_1 \neq e$ then $\pm(d_2 + d_1) = \pm 1$. But since both $d_1$ and $d_2$ are $\geq 1$, this is a contradiction. Therefore $z_1 = e$.

Deduction 2: $||p|| = D \implies ||pt_4|| = D \implies ||a_2 - a_1|| = D \implies (a_1 - a_2)^2 - D = 0$ □

**Theorem 16.** *There is no rank 6 NI-NSD-MTC associated with $G = \langle(01)(2435)\rangle$*

*Proof* If $\sigma = (01)(2435)$, observe, $S_{02} = \epsilon_1\epsilon_2 S_{15} = \epsilon_1\epsilon_2\epsilon_1\epsilon_4 S_{03} \implies \epsilon_0\epsilon_1\epsilon_2\epsilon_4 = 1$.
1) $\epsilon_0 = \epsilon_1 = \epsilon_2 = \epsilon_4$  2) $\epsilon_0 = \epsilon_1 \neq \epsilon_2 = \epsilon_4$  3) $\epsilon_0 = \epsilon_2 \neq \epsilon_1 = \epsilon_4$  4) $\epsilon_0 = \epsilon_4 \neq \epsilon_1 = \epsilon_2$

$$S = \begin{pmatrix} 1 & d_1 & d_2 & d_2 & d_4 & d_4 \\ d_1 & \epsilon_0\epsilon_1 & \epsilon_1\epsilon_4 d_4 & \epsilon_1\epsilon_4 d_4 & \epsilon_1\epsilon_2 d_2 & \epsilon_1\epsilon_2 d_2 \\ d_2 & \epsilon_1\epsilon_4 d_4 & a_1 & a_2 & c_1 & c_2 \\ d_2 & \epsilon_1\epsilon_4 d_4 & a_2 & a_1 & c_2 & c_1 \\ d_4 & \epsilon_1\epsilon_2 d_2 & c_1 & c_2 & \epsilon_2\epsilon_4 a_2 & \epsilon_2\epsilon_4 a_1 \\ d_4 & \epsilon_1\epsilon_2 d_2 & c_2 & c_1 & \epsilon_2\epsilon_4 a_1 & \epsilon_2\epsilon_4 a_2 \end{pmatrix}$$

Case 1 fails as it often will, $0 = S_0 \cdot S_1$ is a sum of positive real numbers ($S_i$ is the $i$th column of $S$). Since columns 2 through 5 have complex numbers we know that column 0 or column 1 is the FPcol. We can always assume the FPcol is another column in the same orbit under the action of $G$, i.e. we can assume that the FPcol is column 0. In case 2, $S_0 \cdot S_1$ implies, $\frac{1}{2} = \frac{d_2 d_4}{d_1}$. But this is a contradiction because $\frac{d_2 d_4}{d_1}$ is an integer based on the $s$-polynomial of $d_2$. Then case 4 is a relabeling of case 3.

After one run of the GBA, we get $c_1 - c_2 = 0$. The action of $G$ on the eigenvalue matrix says that $\frac{c_1}{d_2} \mapsto \frac{-a_2}{d_4} \mapsto \frac{c_2}{d_2} \mapsto \frac{-a_1}{d_4}$. Since $c_1 = c_2 \implies a_1 = a_2$. But then, $c_1 = c_2 = \bar{c}_1$ and $a_1 = a_2 = \bar{a}_1 \implies S$ is real. All NSD-MTC have some non-real entries in $S$. $\square$

**Theorem 17.** *There is no rank 6 NI-NSD-MTC associated with $G = \langle (01)(23), (23)(45) \rangle$*

*Proof* If $\sigma = (01)(23)$, then $\epsilon_0 = \epsilon_1$. This leaves 4 sign choices.
1) $\epsilon_0 = \epsilon_2 = \epsilon_4$  2) $\epsilon_0 = \epsilon_2 \neq \epsilon_4$  3) $\epsilon_0 = \epsilon_4 \neq \epsilon_2$  4) $\epsilon_0 \neq \epsilon_2 = \epsilon_4$

$$S = \begin{pmatrix} 1 & d_1 & d_2 & d_2 & d_4 & d_4 \\ d_1 & 1 & \epsilon_0\epsilon_2 d_2 & \epsilon_0\epsilon_2 d_2 & \epsilon_0\epsilon_4 d_4 & \epsilon_0\epsilon_4 d_4 \\ d_2 & \epsilon_0\epsilon_2 d_2 & a_1 & a_2 & a_5 & \epsilon_2\epsilon_4 a_5 \\ d_2 & \epsilon_0\epsilon_2 d_2 & a_2 & a_1 & \epsilon_2\epsilon_4 a_5 & a_5 \\ d_4 & \epsilon_0\epsilon_4 d_4 & a_5 & \epsilon_2\epsilon_4 a_5 & a_3 & a_4 \\ d_4 & \epsilon_0\epsilon_4 d_4 & \epsilon_2\epsilon_4 a_5 & a_5 & a_4 & a_3 \end{pmatrix}$$

Suppose sign choice 1 or 4 is true. From the orthogonality of the first two columns we can see that $\pm 1 = \frac{d_2^2}{d_1} + \frac{d_4^2}{d_1}$. But by the $s$-polynomials of $d_2$ and $d_4$ we know that $\frac{d_2^2}{d_1}$ and $\frac{d_4^2}{d_1}$ are integers. But 1 is not a sum of two positive or two negative integers. Then, sign choice 3 is a relabeling of sign choice 2.

Note that if $a_5 = 0$, then the $S$-matrix is identical to $S1$ in theorem 13. That theorem was based on theorem 15. The work for theorem 10 doesn't require $(01) \in G$, but that the orbit of 0 is $\{0, 1\}$ and the existence of the appropriate 0's. And, $a_5 = 0$ gives us the needed 0's for the results of theorem 10 to still be true. So, the work done eliminating $S1$ holds here. Thus, we can assume $a_5 \neq 0$. Table 2 gives the details of each GBA calculation for sign choice 2.

Deduction: In the output of $(01)(23),(23)(45)$-2-4, many of the relevant polynomials had $p - 1$ as a factor. It's always true that $||p||^2 = D$. Therefore, $p - 1 = 0 \implies D - 1 = 0$. But, since only columns 0 and 1 have all real entries, the FPcol is either column 0 or 1. Since both are in the same orbit, we may assume the FPcol is column 0 and therefore that $d_i \geq 1 \implies D \geq 6$. Thus, $p - 1 \neq 0$. If $t_4 + 1 \neq 0 \implies t_2 + 1 = t_4^2 - t_4 + 1 = 0 \implies N = 6$. But, $G$ cannot be a subgroup of $\text{Gal}(\mathbb{Q}_N/\mathbb{Q})$ for $N = 2$ or $N = 6$, as neither have a subgroup of order 4. $\square$

**Theorem 18.** *There is no rank 6 NI-NSD-MTC associated with $G = \langle (0213), (45) \rangle$.*

*Proof* Let $\sigma = (0213)$, then $S_{04} = \epsilon_2\epsilon_4 S_{24} = \epsilon_1\epsilon_2 S_{14} = \epsilon_1\epsilon_2\epsilon_3\epsilon_4 S_{34} = \epsilon_0\epsilon_1\epsilon_2\epsilon_3 S_{04} \implies \epsilon_0\epsilon_1\epsilon_2\epsilon_3 = 1$. This gives 8 sign choices.

1) $\epsilon_0 = \epsilon_1 = \epsilon_2 = \epsilon_3 = \epsilon_4$  2) $\epsilon_0 = \epsilon_1 = \epsilon_2 = \epsilon_3 \neq \epsilon_4$  3) $\epsilon_0 = \epsilon_1 \neq \epsilon_2 = \epsilon_3 = \epsilon_4$
4) $\epsilon_0 = \epsilon_1 = \epsilon_4 \neq \epsilon_2 = \epsilon_3$  5) $\epsilon_0 = \epsilon_2 \neq \epsilon_1 = \epsilon_3 = \epsilon_4$  6) $\epsilon_0 = \epsilon_2 = \epsilon_4 \neq \epsilon_1 = \epsilon_3$
7) $\epsilon_0 = \epsilon_3 \neq \epsilon_1 = \epsilon_2 = \epsilon_4$  8) $\epsilon_0 = \epsilon_3 = \epsilon_4 \neq \epsilon_1 = \epsilon_2$

Sign choices 1 and 2 give at least 2 columns with all positive entries. Choice 4 is a relabeling of choice 3. Similarly, choice 8 is a relabeling of choice 7.

$$S = \begin{pmatrix} 1 & d_1 & d_2 & d_3 & d_4 & d_4 \\ d_1 & 1 & \epsilon_0\epsilon_1 d_3 & \epsilon_0\epsilon_1 d_2 & \epsilon_1\epsilon_2 d_4 & \epsilon_1\epsilon_2 d_4 \\ d_2 & \epsilon_0\epsilon_1 d_3 & \epsilon_0\epsilon_3 d_1 & \epsilon_0\epsilon_2 & \epsilon_2\epsilon_4 d_4 & \epsilon_2\epsilon_4 d_4 \\ d_3 & \epsilon_0\epsilon_1 d_2 & \epsilon_0\epsilon_2 & \epsilon_0\epsilon_3 d_1 & \epsilon_0\epsilon_4 d_4 & \epsilon_0\epsilon_4 d_4 \\ d_4 & \epsilon_1\epsilon_2 d_4 & \epsilon_2\epsilon_4 d_4 & \epsilon_0\epsilon_4 d_4 & a_1 & a_2 \\ d_4 & \epsilon_1\epsilon_2 d_4 & \epsilon_2\epsilon_4 d_4 & \epsilon_0\epsilon_4 d_4 & a_2 & a_1 \end{pmatrix}$$

The non-self-dual columns must contain some non-real entry. The self dual columns must contain only real entries. Therefore, $a_1$ and $a_2$ are non-real. Tables 3, 4, 5, and 6 give the details of each GBA calculation for sign choices 3, 5, 6, and 7 respectively. □

**Theorem 19.** *There is no rank 6 NI-NSD-MTC associated with $G = \langle (01)(23), (02)(13), (45) \rangle$.*

*Proof* Let $\sigma = (01)(23)$ and $\tau = (02)(13)$. Then $\epsilon_0 = \epsilon_1, \epsilon_2 = \epsilon_3, \epsilon_4 = \epsilon_5$ and $\delta_0 = \delta_2, \delta_1 = \delta_3, \delta_4 = \delta_5$. Note that

$$S_{03} = \epsilon_0\epsilon_2 S_{12} = \delta_0\delta_1 S_{21} \implies \epsilon_0\epsilon_2 = \delta_0\delta_1$$

.

$$S = \begin{pmatrix} 1 & d_1 & d_2 & d_3 & d_4 & d_4 \\ d_1 & 1 & \epsilon_0\epsilon_2 d_3 & \epsilon_0\epsilon_2 d_2 & \epsilon_0\epsilon_4 d_4 & \epsilon_0\epsilon_4 d_4 \\ d_2 & \epsilon_0\epsilon_2 d_3 & 1 & \epsilon_0\epsilon_2 & \delta_0\delta_4 d_4 & \delta_0\delta_4 d_4 \\ d_3 & \epsilon_0\epsilon_2 d_2 & \epsilon_0\epsilon_2 d_1 & 1 & \epsilon_0\epsilon_4\delta_1\delta_4 d_4 & \epsilon_0\epsilon_4\delta_1\delta_4 d_4 \\ d_4 & \epsilon_0\epsilon_4 d_4 & \delta_0\delta_4 d_4 & \epsilon_0\epsilon_4\delta_1\delta_4 d_4 & a_1 & a_2 \\ d_4 & \epsilon_0\epsilon_4 d_4 & \delta_0\delta_4 d_4 & \epsilon_0\epsilon_4\delta_1\delta_4 d_4 & a_2 & a_1 \end{pmatrix}$$

This leaves 8 cases,
1) $\epsilon_0 = \epsilon_2 = \epsilon_4$    2) $\epsilon_0 = \epsilon_2 = \epsilon_4$    3) $\epsilon_0 = \epsilon_2 \neq \epsilon_4$    4) $\epsilon_0 = \epsilon_2 \neq \epsilon_4$
   $\delta_0 = \delta_1 = \delta_4$       $\delta_0 = \delta_1 \neq \delta_4$       $\delta_0 = \delta_1 \neq \delta_4$       $\delta_0 = \delta_1 = \delta_4$
5) $\epsilon_0 = \epsilon_4 \neq \epsilon_2$  6) $\epsilon_0 = \epsilon_4 \neq \epsilon_2$  7) $\epsilon_0 \neq \epsilon_2 = \epsilon_4$  8) $\epsilon_0 \neq \epsilon_2 = \epsilon_4$
   $\delta_0 = \delta_4 \neq \delta_1$       $\delta_0 \neq \delta_1 = \delta_4$       $\delta_0 \neq \delta_1 = \delta_4$       $\delta_0 = \delta_4 \neq \delta_1$

Since columns 4 and 5 must have non-real entries, the FPcol must be either 0,1,2 or 3. But they share an orbit, so we may assume that the FPcol is column 0 and that $d_i > 0$ for all $i$. Cases 1-4 leave at least two columns with all positive entries, but this is only allowable for the FPcol. Cases 6 and 8 are relabelings of case 5.

Note, if $d_2 - 1 = 0$, then $d_1 + d_3 = 0$, because $\frac{-d_3}{d_1}$ is a Galois conjugate of $d_2$. But we can assume that all $d_i \geq 1$. Therefore $d_2 \neq 1$. Then, $h_2 = -D + (d_3 + 1)^2$ and $h_3 = -D + (d_2 - 1)^2 \implies (d_3 + 1)^2 = (d_2 - 1)^2$. But $d_i geq 1 \implies d_3 + 1 = d_2 - 1$. Then $h_3 = d_3 + 1 - d_2 + 1 = 0$. In sign choice 7, $h_1 + h_3 = d_1 - d_3 = 0$. But $\frac{-d_3}{d_1} = -1$ is a Galois conjugate of $d_2$. Tables 7 and 8 gives the details of each GBA calculation for sign choices 5 and 7 respectively. □

**Theorem 20.** *There is no rank 6 NI-NSD-MTC associated with $G = \langle (012), (45) \rangle$.*

*Proof* We may assume that $\epsilon_0 = \epsilon_2 \neq \epsilon_1 = \epsilon_4$. So,

$$S = \begin{pmatrix} 1 & d_1 & d_2 & d_3 & d_4 & d_4 \\ d_1 & -d_2 & -1 & d_3 & d_4 & d_4 \\ d_2 & -1 & d_1 & -d_3 & -d_4 & -d_4 \\ d_3 & d_3 & -d_3 & z_1 d_3 & z_2 d_4 & z_2 d_4 \\ d_4 & d_4 & -d_4 & z_3 d_3 & a_1 & a_2 \\ d_4 & d_4 & -d_4 & z_3 d_3 & a_2 & a_1 \end{pmatrix}$$

One run of the GBA yields, $t_3 - t_4 = 0$. Adding that relation and running the GBA one more time yields, $d_4 z_3 - a_1 = 0$. But that forces $a_1$ to be real and it must be non-real. $\square$.

**Theorem 21.** *There is no rank 6 NI-NSD-MTC associated with $G = \langle (01)(23), (45) \rangle$.*

*Proof* Let $\sigma = (01)(23)$. Since $\sigma^2 = Id$ and $\sigma$ fixes labels 4 and 5, $\epsilon_0 = \epsilon_1$ and $\epsilon_2 = \epsilon_3$. Since $4^* = 5 \implies \epsilon_4 = \epsilon_5$. This gives 4 sign choices,
1) $\epsilon_0 = \epsilon_2 = \epsilon_4$  2) $\epsilon_0 = \epsilon_2 \neq \epsilon_4$  3) $\epsilon_0 \neq \epsilon_2 = \epsilon_4$  4) $\epsilon_0 = \epsilon_4 \neq \epsilon_2$

$$S = \begin{pmatrix} 1 & d_1 & d_2 & d_3 & d_4 & d_4 \\ d_1 & 1 & \epsilon_0\epsilon_2 d_3 & \epsilon_0\epsilon_2 d_2 & \epsilon_0\epsilon_4 d_4 & \epsilon_0\epsilon_4 d_4 \\ d_2 & \epsilon_0\epsilon_2 d_3 & a_1 & a_2 & a_3 & a_3 \\ d_3 & \epsilon_0\epsilon_2 d_2 & a_2 & a_1 & \epsilon_2\epsilon_4 a_3 & \epsilon_2\epsilon_4 a_3 \\ d_4 & \epsilon_0\epsilon_4 d_4 & a_3 & \epsilon_2\epsilon_4 a_3 & a_4 & a_5 \\ d_4 & \epsilon_0\epsilon_4 d_4 & a_3 & \epsilon_2\epsilon_4 a_3 & a_5 & a_4 \end{pmatrix}$$

In case 1, there will always be at least two columns that are all positive or all negative. We cannot assume that the FPcol is column 0. However, we can assume it's either column 0 or column 2. If it's column 0, then both column and 0 and column 1 are all positive and $0 = S_0 \cdot S_1$ is therefore a sum of positive real numbers. If the FPcol is column 2, then column 2 is entirely positive or negative. But column three has the same entries but in a different order. In either choice ($\pm$), $0 = S_2 \cdot S_3$ is a sum of positive real numbers.

In all three remaining cases we'll show that $1, t_1, t_2, t_3,$ and $t_4$ are all distinct, and that $t_1 + 1$ and $t_2 + t_3$ are not both 0 at the same time. This is enough, but we'll need some results from the next section to say why.

In all cases we suppose $t_1 - 1 = 0$. This is the most difficult case to eliminate. Progressing through successive runs of the GBA often yields products such as, $(t_2 - t_3)(d_2 z_2 - d_3 z_3) = 0$ ($z_2$ and $z_3$ are the coefficients of the s-polynomial of $d_3$). Supposing that $t_2 - t_3 = 0$ eventually yields that the category must be integral. Therefore we assume $d_2 z_2 - d_3 z_3 = 0$. In sign choice 3 such products were encountered. Sign choices 3 and 4 had 2 each. After each of those it was much more straight forward. If $t_2 - 1 = t_3 - t_1 = 0$, then in three runs of the GBA we get $pD = 0$. If $t_2 - t_3 = 0$ we get a contradiction in two runs of the GBA. Since we're free to switch the labels 2 and 3, This shows that $1, t_1, t_2, t_3, t_4$ are distinct. Next suppose $t_1 + 1 = t_2 + t_3 = 0$. Then, the GBA immediately yields that $p = 0$. This is indeed enough to show there are no NI-NSD-MTC of rank 6 with Galois group $\langle (01)(23), (45) \rangle$, however we need the results from lemmas 25 and 26 to prove it.

Table 1. (01)

| Step: | Zero Factor |
|---|---|
| S1-1 | N/A |
| Contradiction: | Empty Variety |
| S2-1 | $pD$ |
| Contradiction: | $p, D \neq 0$ |
| S3-1 | $p^2 + D$ |
| | $d_2 z_1 - d_1 e + 1$ |
| | $z_1^2 - 1$ |
| Deduction 1 | $z_1 - e = 0$ |
| S3-2 | $2pt_2 - 3d_2 e + p$ |
| | $2d_2 t_2 + pe + d_2$ |
| | $d_3 z_2 + 2d_4 z_4 + d_2 e + 2$ |
| S3-3 | $pt_4^3 - a_1 + a_2$ |
| Deduction 2 | $(a_1 - a_2)^2 - D = 0$ |
| S3-4 | $D$ |
| Contradiction: | $D \neq 0$ |

Table 2. (01)(23),(23)(45)-2

| Gröbner Basis Run | Zero Factor |
|---|---|
| (01)(23),(23)(45)-2-1 | $a_1 - a_2 + a_3 - a_4$ |
| | $d_1 - a_3 - a_4 - 1$ |
| | $d_1 + a_1 + a_2 + 1$ |
| | $t_1 - 1$ |
| (01)(23),(23)(45)-2-2 | $p^2 - D$ |
| | $6a_5^2 + D$ |
| (01)(23),(23)(45)-2-3 | $D - 5p + 4$ |
| (01)(23),(23)(45)-2-4 | $p - 4$ |
| | $(t_2 + 1)(t_4 + 1)$ |
| | $(t_4^2 - t_4 + 1)(t_4 + 1)$ |
| Deduction | $t_4 + 1 = 0$ |
| (01)(23),(23)(45)-2-5 | $d_2$ |
| Contradiction: | $d_2 \neq 0$ |

Table 3. (0213),(45)-3

| Gröbner Basis Run | Zero Factor(s) |
|---|---|
| (0213),(45)-3-1 | $h_1$ |
| (0213),(45)-3-2 | $h_2, h_3, h_4$ |
| (0213),(45)-3-3 | $h_5, h_6$ |
| (0213),(45)-3-4 | $h_7, h_8, h_9$ |
| (0213),(45)-3-5 | $d_2 - d_3 + a_2$ |
| Contradiction: | $a_2$ is non-real |

Table 4. (0213),(45)-5

| Gröbner Basis Run | Zero Factor(s) |
|---|---|
| (0213),(45)-5-1 | $h_1, \ldots, h_8$ |
| (0213),(45)-5-2 | $h_9$ |
| (0213),(45)-5-3 | $d_2 + d_3 - a_2$ |
| Contradiction: | $a_2$ is non-real |

Table 5. (0213),(45)-6

| Gröbner Basis Run | Zero Factor(s) |
|---|---|
| (0213),(45)-6-1 | $h_1, \ldots, h_{10}$ |
| (0213),(45)-6-2 | $h_{11}$ |
| (0213),(45)-6-3 | $d_2 + d_3 + a_2$ |
| Contradiction: | $a_2$ is real |

Table 6. (0213),(45)-7

| Gröbner Basis Run | Zero Factor(s) |
|---|---|
| (0213),(45)-7-1 | $h_1, \ldots, h_8$ |
| (0213),(45)-7-2 | $h_9$ |
| (0213),(45)-7-3 | $h_{10}, h_{11}, h_{12}$ |
| (0213),(45)-7-4 | $d_2 - d_3 + a_2$ |
| Contradiction: | $a_2$ is non-real |

Table 7. (01)(23),(02)(13),(45)-5

| Gröbner Basis Run | Zero Factor(s) |
|---|---|
| (01)(23),(02)(13),(45)-5-1 | $h_1, h_2$ |
|  | $(d_2 - 1)(-D + (d_2 - 1)^2)$ |
| Deduction | $h_3 = 0$ |
| (01)(23),(02)(13),(45)-5-2 | $D$ |
| Contradiction: | $D \neq 0$ |

Table 8. (01)(23),(02)(13),(45)-7

| Gröbner Basis Run | Zero Factor (s) |
|---|---|
| (01)(23),(02)(13),(45)-7-1 | $h_1, h_2$ |
| (01)(23),(02)(13),(45)-7-2 | $h_3$ |
| Deduction | $d_2 + 1 = 0$ |
| Contradiction: | $d_2 > 0$ |

# 5. Non-Integral, Self-Dual, Modular Categories

## 5.1. Representation Theory.
Often the orthogonality and twist relations aren't enough. But there is a way to deduce a set of possible relations on the entries of the $T$ matrix through modular representations. As an abstract group, $\mathrm{SL}(2,\mathbb{Z}) \cong \langle \mathfrak{s}, \mathfrak{t} \mid \mathfrak{s}^4 = 1, \ (\mathfrak{s}\mathfrak{t})^3 = \mathfrak{s}^2 \rangle$. The standard choice is,

$$\mathfrak{s} := \begin{bmatrix} 0 & -1 \\ 1 & 0 \end{bmatrix} \text{ and } \mathfrak{t} := \begin{bmatrix} 1 & 1 \\ 0 & 1 \end{bmatrix}$$

If $(S, T)$ are the modular data of a category, $\mathcal{C}$, then let $\eta : \mathrm{GL}(\Pi_\mathcal{C}, \mathbb{C}) \to \mathrm{PGL}(\Pi_\mathcal{C}, \mathbb{C})$ be the natural surjection. Then, $\bar\rho(\mathfrak{s}) = \eta(S)$ and $\bar\rho(\mathfrak{t}) = \eta(T)$ defines a projective representation of $\mathrm{SL}(2,\mathbb{Z})$.

**Definition 6.** *A **modular representation**, $\rho : GL(\Pi_\mathcal{C}, \mathbb{C}) \to GL(\Pi_\mathcal{C}, \mathbb{C})$, is a representation of $SL(2,\mathbb{Z})$ such that, $\bar\rho = \eta \circ \rho$. We call $(s,t)$, where $s = \rho(\mathfrak{s})$ and $t = \rho(\mathfrak{t})$ a **normalized modular pair**. [3]*

Modular representations do exist, [5] gives a construction. Moreover, they construct the complete set. In any modular representation, $\mathfrak{t} \mapsto \frac{x}{\zeta}T$, where $x^{12} = 1$ and $\zeta^6 = \frac{p_+}{p_-} \implies \frac{x}{\zeta}$ is a root of unity. Therefore, $t = \rho(\mathfrak{t}) = \gamma T$ for some root of unity $\gamma$.

**Theorem 22.** *Let $\mathcal{C}$ be a modular category of rank $r$, with $T$-matrix of order $N$. Suppose $(s,t)$ is normalized modular pair of $\mathcal{C}$. Set $t = (\delta_{ij}t_i)$ and $n = \mathrm{ord}(t)$. Then,*
*(a) $N \mid n \mid 12N$ and $s, t \in Gl_r(\mathbb{Q}_n)$. Moreover,*
*(b) (Galois Symmetry) for $\sigma \in Gal(\mathbb{Q}_n/\mathbb{Q})$, $\sigma^2(t_i) = t_{\sigma(i)}$ [3]*

Here we begin to see where we can get information about $T$. The next theorem will tell us why, but it's very useful to simply know which entries of $T$ are distinct. And since $t = \gamma T$, if the entries of $T$ are distinct then so are the entries of $t$.

**Theorem 23.** *Let $\mathcal{C}$ be a modular category of rank, $r$, and $\rho : SL(2,\mathbb{Z}) \to GL(r, \mathbb{C})$ a modular representation of $\mathcal{C}$. Then $\rho$ cannot be isomorphic to a direct sum of two representations with disjoint t-spectra. [3]*

Now we see that whether the entries of $T$ overlap affects the structure of a modular representation. But the entries of $T$ don't just affect the structure of the representation in that way, the order of the $T$ matrix also has an important impact.

**Theorem 24.** *Let $(S, T)$ be the modular data of the modular category $\mathcal{C}$ with $N = \mathrm{ord}(T)$. Then $N$ is minimal such that the projective representation $\bar\rho_\mathcal{C}$ of $SL(2,\mathbb{Z})$ associated with the modular data can be factored through $SL(2, \mathbb{Z}/N\mathbb{Z})$. [3]*

A lot is known about $\mathrm{SL}(2, \mathbb{Z}/p^\lambda\mathbb{Z})$ representations in low rank. In fact usually the $t$-spectra is known for all representations when $p$ and $r$ and small enough. The $t$-spectra may be known, but we cannot make any immediate assumptions about the ordering, i.e. we will not know which element of the spectra is $t_0$. Note that since $T_0 = 1$, then clearly $t_0 = \gamma$. Fortunately, while we may not be able to deduce exactly what any particular entry of $T$ is, we can deduce some relations on $T$, that may or may not depend on which element is $\gamma$.

Sometimes, we will not directly appeal to a particular representation. Eholzer who has compiled most of the relevant $t$-spectra in various tables in [6] and [7]. Bruillard, Ng, Rowell, and Wang have taken parts of Eholzer's tables and condensed them into a very helpful and more explicit list in [3]. Most helpful though is Eholzer's table 12 in [7]. In it are the all simple and non-degenerate strong fusion algebras. For this paper this means that if a group has completely disjoin $t$-spectra (this is the non-degeneracy condition) and also that we can show it isn't the Galois group of a product category (this is the simple condition) that's enough.

To apply these theorems, we first determine how the $T$ spectra may overlap based on Galois symmetry. For example, if $\sigma = (01)(23)$ and $\theta_0 = \theta_2$, then Galois symmetry says that $\theta_1 = \theta_3$. We try to eliminate cases with the most overlap first. Usually there's a very small number of them and one or two runs of the GBA is enough to find a contradiction. In most cases we will be able to conclude there is exactly $k$ pairs of $i,j$ such that $\theta_i = \theta_j$ and $i \neq j$, and then draw strong conclusions about the number of possible direct summands of $\rho$. Recall that $G$ is a subgroup of $\text{Gal}(\mathbb{Q}_N/\mathbb{Q})$ and that $\text{Gal}(\mathbb{Q}_N/\mathbb{Q})/G \cong (\mathbb{Z}/2\mathbb{Z})^k$ for some integer $k$. This severely restricts the possible values of $N$ and $n$. Then we must use the known representations to attempt to piece together the possible $t$-spectras of $\rho$. For shorthand, we'll use the dimension to represent the subrepresentations. For example a 6 dimensional representation may break down as a direct sum of two subrepresentations in the following ways, $1 \oplus 5$, $2 \oplus 4$, and $3 \oplus 3$. But remember the subrepresentations must be factored to get to the $\mathbb{Z}/p^{\lambda}\mathbb{Z}$ irreducible representations. This means that $2 \oplus 4$ might break down as a $2 \oplus 2 \otimes 2$. If we find a combination of irreducible representations that meet all the criteria, we then deduce whatever relations on the $t$-spectra we can and add them the GBA. If we do not find a combination, then we have the contradiction we need to eliminate the group.

**Lemma 25.** *If $n$ divides 48, and $\sigma \in \text{Gal}(\mathbb{Q}_n/\mathbb{Q})$ such that $\sigma$ has order 4. Then for any primitive 16th or 48th root of unity, $\zeta \in \mathbb{Q}_n$, $\sigma^2(\zeta) = \pm \zeta$.*

*Proof:* Let $\sigma, \tau \in \text{Gal}(\mathbb{Q}_n, \mathbb{Q}) \cong (\mathbb{Z}/2\mathbb{Z})^{k_n} \oplus \mathbb{Z}/4\mathbb{Z}$, where $k_{16} = 1$ and $k_{48} = 2$. Suppose both $\sigma$ and $\tau$ are of order 4, then when viewed under the isomorphism, it's clear that $\sigma^2 = \tau^2$. To prove the lemma, we need only prove it's true for one element of order 4. Moreover, suppose that $\sigma^2(\zeta_n) = -\zeta_n$, where $\zeta_n = e^{2\pi i/n}$. Then note that all other primitive $n^{th}$ roots of unity are odd powers of $\zeta_n$. Therefore, $\sigma^2(\zeta_n^{2l+1}) = (-1)^{2l+1} \zeta_n^{2l+1} = -\zeta_n^{2l+1}$. So, we need only show that for one $\sigma$, $\sigma^2(\zeta_n) = -\zeta_n$. If $n = 16$, then let $\sigma(\zeta_{16}) = \zeta_{16}^3 \implies \sigma^2(\zeta_{16}) = \zeta_{16}^9 = -\zeta_{16}$. Such a $\sigma$ exists because $\zeta_{16}^3$ is also a primitive $16^{th}$ root of unity. Similarly, if $n = 48$, then let $\sigma(\zeta_{48}) = \zeta_{48}^5 \implies \sigma^2(\zeta_{48}) = \zeta_{48}^{25} = -\zeta_{48}$. $\square$

**Lemma 26.** *If 5 divides $n$ and 16 does not and there are at most two pair of identical $t_i$, then*
  *1) The reduced and factored of any modular representation must fit the following dimensions, $\rho = 6$, $1 \oplus 5$, $2 \oplus 4$, $2 \oplus 2 \otimes 2$, $1 \oplus 1 \oplus 4$, $1 \oplus 1 \oplus 2 \otimes 2$, or $3 \oplus 3$.*
  *2) If $\rho = 6$ or $1 \oplus 5$, then the $t$-spectra $= \alpha \otimes \{1, 1, \zeta_5, \zeta_5^4, \zeta_5^2, \zeta_5^3\}$*
  *3) If $\rho = 2 \oplus 4$ or $1 \oplus 1 \oplus 4$, then the $t$-spectra $= \alpha \otimes \{\zeta_5, \zeta_5^4, \zeta_5, \zeta_5^4, \zeta_5^2, \zeta_5^3\}$*
  *4) If $\rho = 2 \oplus 2 \otimes 2$ or $1 \oplus 1 \oplus 2 \otimes 2$, then the $t$-spectra $=$*
$\alpha \otimes \{\alpha_1 \zeta_5, \alpha_1 \zeta_5^4, \alpha_1 \zeta_5, \alpha_1 \zeta_5^4, \alpha_2 \zeta_5, \alpha_2 \zeta_5^4\}$.
  *5) If $\rho = 3 \oplus 3$, then the $t$-spectra is the same as 1) or the $t$-spectra $= \alpha \otimes \{1, 1, \alpha_1, \alpha_2, \zeta_5, \zeta_5^4\}$*
*Where, $\zeta_5$ is any primitive $5^{th}$ root of unity, $\alpha$ is some unknown $24^{th}$ root of unity, and $\alpha_1, \alpha_2$ are both a, $2^{nd}$, $3^{rd}$, $4^{th}$, or $8^{th}$ root of unity.*

*Proof:* The first point follows immediately from theorem 5.3. There are at most 3 possible summands of $\rho$ when reduced. Otherwise, we could group them in such away to have $\rho$ as the direct summand of two subrepresentations without overlapping $t$-spectra. As a reminder, $2 \otimes 2$ is a 4 dimensional representation, that has been factored into two different prime powers. In all the cases listed in point 1, there could be other 1 dimensional factors. Those have been left off for simplicity, as the restriction on overlapping $t$-spectra from theorem 5.3 almost always forces them to be removed when factoring $\gamma$ out of the $t$-spectra to get the $T$-spectra.

To work on the remaining points, we need to note that if $\zeta_5$ appears anywhere then a corresponding $\zeta_5^4$ must appear. This is true because of Galois symmetry. Let $\sigma \in \text{Gal}(\mathbb{Q}_n/\mathbb{Q})$ such that $\sigma$ has order 4. Then if $\zeta_5$ appears $\sigma^2(\zeta_5) = \zeta_5^4$ must appear in a similar manner, i.e. if $\alpha \zeta_5$ is in the $t$-spectra, then $\alpha \zeta_5^4$ must be too. For simplicity we'll use $\zeta_5$ and $\zeta_5^4$ unless we need to distinguish two different pairs of $5^{th}$ roots. Finally, since 5 divides $n$, one such pair exist, i.e. we only consider the representations of dimension greater than 1. Then by [6], the possible order representations are, $\{\zeta_5, \zeta_5^4\}$, $\{1, \zeta_5, \zeta_5^4\}$, $\{\zeta_5, \zeta_5^4 \zeta_5^2, \zeta_5^3\}$, $\{1, \zeta_5, \zeta_5^4, \zeta_5^2, \zeta_5^3\}$, and $\{1, 1, \zeta_5, \zeta_5^4, \zeta_5^2, \zeta_5^3\}$.

The points then follow from combining these representations in all allowable ways. Note that for $2 \otimes 2$, this must break down as $\{\alpha_1, \alpha_2\} \otimes \{\zeta_5, \zeta_5^4\}$, where $\alpha_1$ and $\alpha_2$ are both prime power roots of unity. For example they couldn't be primitive $6^{th}$ roots of unity. But $\alpha$ could be a primitive $6^{th}$ root since it could be a product of multiple different 1 dimensional factors. However, that will have no affect on the $T$-matrix as it will be factored out with $\gamma$.

It may appear that we left off one possibility, $2 \oplus 2 \oplus 2$. A priori, this is a possible decomposition of $\rho$ with the assumed overlap in the $t$-spectra. However, $\zeta_5$ and $\zeta_5^4$ are forced to appear in all three summands. This forces the $t$-spectra to have 3 identical pairs of $\gamma t_i$'s. $\square$

## 5.2. Non-Self-Dual. Proof of Theorem 19 continued:

We can now finish off the proof of theorem 19. We know that $1, t_1, t_2, t_3, t_4$ are distinct. Galois symmetry tells us that we have some element of order 4 in $\text{Gal}(\mathbb{Q}_n/\mathbb{Q})$ and that's the exponent of the group. Cyclotomic extensions are well known. The only cyclotomic Galois groups with elements of order 4 and no higher, means either 5 or 16 divides $n$. We also know that all elements of order 4 are liftings of $\sigma = (01)(23)$. If both 5 and 16 divide $n$ the corresponding elements of order 4 would act differently on the $t$ spectra. This forces only 5 or 16 not both to divide $n$. If 16 divides $n$, then by lemma 25 and Galois symmetry, $t_1 + 1 = t_2 + t_3 = 0$. We've already eliminated that case. Looking at lemma 26, specifically the ones with at most 1 pair of identical $t_i$, we see that either $N = 5$ or the $t$-spectra is $\{1, 1, \alpha_1, \alpha_2, \zeta_5, \zeta_5^4\}$. Suppose $N = 5$. Observe that $\text{Gal}(\mathbb{Q}_N/\mathbb{Q})$ has no subgroup isomorphic to $\mathbb{Z}/2\mathbb{Z} \times \mathbb{Z}/2\mathbb{Z}$. Then suppose the $t$-spectra is $\{1, 1, \alpha_1, \alpha_2, \zeta_5, \zeta_5^4\}$. Let $\Sigma$ be a lifting of $\sigma = (01)(23)$. Since $t_1 \neq 1$ and $t_2 \neq t_3$, $\Sigma^2$ moves $\gamma$, $\gamma t_1$, $\gamma t_2$, and $\gamma t_3$. But, all elements squared will fix, 1, 1, $\alpha_1$, and $\alpha_2$. Therefore the Galois symmetry doesn't lineup. And 5 can't divide $n$ either. $\square$

**Theorem 27.** *All rank 6 non-integral, non-self-dual categories are product categories.*

*Proof* Already shown, all rank 6 NI-NSD-MTC, have $\langle (01)(24)(35), (23)(45) \rangle$ as their Galois group up to a relabeling of the classes of simple objects. Let $\sigma = (01)(24)(35) \rangle$.

$$S = \begin{pmatrix} 1 & d_1 & d_2 & d_2 & d_4 & d_4 \\ d_1 & \epsilon_0 \epsilon_1 & \epsilon_1 \epsilon_4 d_4 & \epsilon_1 \epsilon_4 d_4 & \epsilon_1 \epsilon_2 d_2 & \epsilon_1 \epsilon_2 d_2 \\ d_2 & \epsilon_1 \epsilon_4 d_4 & a_1 & a_2 & a_3 & a_4 \\ d_2 & \epsilon_1 \epsilon_4 d_4 & a_2 & a_1 & a_4 & a_3 \\ d_4 & \epsilon_1 \epsilon_2 d_2 & a_3 & a_4 & \epsilon_2 \epsilon_4 a_1 & \epsilon_2 \epsilon_4 a_2 \\ d_4 & \epsilon_1 \epsilon_2 d_2 & a_4 & a_3 & \epsilon_2 \epsilon_4 a_2 & \epsilon_2 \epsilon_4 a_1 \end{pmatrix}$$

Observe that $S_{02} = \epsilon_1 \epsilon_2 S_{14} = \epsilon_1 \epsilon_2 \epsilon_0 \epsilon_4 S_{02}$. Therefore, $\epsilon_0 \epsilon_1 \epsilon_2 \epsilon_4 = 1$. Suppose $\epsilon_0 = \epsilon_1 \implies \epsilon_2 = \epsilon_4$. Then check orthogonality of the first two columns to see that, $2d_1 + \epsilon_0 \epsilon_2 4 d_2 d_4 = 0 \implies \frac{-\epsilon_0 \epsilon_2}{2} = \frac{d_2 d_4}{d_1} \in \mathbb{Z}$ by the $s$-polynomial of $d_2$. Thus we get that $\epsilon_0 \neq \epsilon_1$ and $\epsilon_2 \neq \epsilon_4$. Since we can relabel the dual pairs, $2 \mapsto 4$ and $3 \mapsto 5$, we can suppose $\epsilon_4 = \epsilon_0$.

The typical GBA is not immediately helpful. So we need to try to understand the entries of the $T$ matrix in order to apply the representation theory theorems. To that end we need to show that, $1, t_1, t_2,$ and $t_4$ are distinct.

Suppose, $t_1 - 1 = 0$ then one run of the GBA yields $pD = 0$. Let $t_2 - 1 = 0 \implies t_4 - t_1 = 0$. Then the GBA yields $t_4^2 - 1 = 0 \implies N = 1$ or 2 (Galois Symmetry). Let $t_4 - 1 = 0 \implies t_2 - t_1 = 0$. The GBA then yields $t_2^2 - 1 = 0$, again $N = 1$ or 2. Let $t_2 - t_4 = 0$. The GBA yields, $t_1 - 1 = 0$, but we've already shown that doesn't happen. Therefore, 5 or 16 divide $n$ but not both and there at most 2 pair of identical $t_i$. Suppose 16 divides $n$. Then, $t_1 + 1 = t_2 + t_4 = 0$. After two runs of the GBA, we get $t_4^3 - a_1 - a_2 = 0 \implies t_4^3$ s real. Then $N$ is 2 or 6. But then $\text{Gal}(\mathbb{Q}_N\mathbb{Q})$ would not have a subgroup of order 4. Thus, 5 divides $n$.

We use the possible $t$-spectra from lemma 26, $\alpha \otimes \{\zeta_5, \zeta_5^4, \zeta_5, \zeta_5^4, \zeta_5^2, \zeta_5^3\}$ and $\alpha \otimes \{\alpha_1\zeta_5, \alpha_1\zeta_5^4, \alpha_1\zeta_5, \alpha_1\zeta_5^4, \alpha_2\zeta_5, \alpha_2\zeta_5^4\}$. These are the only two that have two pair of identical $t_i$.

Suppose the $t$-spectra is $\alpha \otimes \{\zeta_5, \zeta_5^4, \zeta_5, \zeta_5^4, \zeta_5^2, \zeta_5^3\}$. Then the $T$-spectra is $\{1, \zeta_5, \zeta_5^2, \zeta_5^2, \zeta_5^4, \zeta_5^4\}$. Then, $t_1^4 + t_1^3 + t_1^2 + t_1 + 1 = (t_1^2 - t2)(t_1^2 - t4) = (t_2^2 - t_4)(t_4^2 - t_2) = 0$. In both the products the factors on the left occur at the same time as do the factors on the right. In either case, one run of the GBA yields, $p = 0$.

Suppose the $t$-spectra is $\alpha \otimes \{\alpha_1\zeta_5, \alpha_1\zeta_5^4, \alpha_1\zeta_5, \alpha_1\zeta_5^4, \alpha_2\zeta_5, \alpha_2\zeta_5^4\}$. Then the $T$-spectra is $\{1, \zeta_5^3, \beta, \beta, \beta\zeta_5^3, \beta\zeta_5^3\}$, where $\beta = \alpha_1/\alpha_2$. Then, $t_1^4 + t_1^3 + t_1^2 + t_1 + 1 = (t_2t_1 - t_4)(t_4t_1 - t_2) = 0$. Note that $\alpha_1$ and $\alpha_2$ came from the same prime power representation. Therefore $\beta$ is a primitive $2^{nd}$, $3^{rd}$, $4^{th}$ or $8^{th}$ root of unity.

Let $t_4t_1 - t_2 = 0$. Since $t_4 = \beta$, $(t_4 + 1)(t_4^2 + 1)(t_4^2 + t_4 + 1)(t_4^4 + 1) = 0$. If $t_4 + 1 = 0$ or $t_4^2 + 1 = 0$, one run of the GBA yields $p = 0$. If $t_4^2 + t_4 + 1 = 0$, two runs of the GBA yield, $d_4 - 1 = 0$. But then $d_2 = -d_1$. Since columns 2 through 5 necessarily have some non-real entries, we can assume the FPcol is column 0 and that $d_i > 0$ for all $i$. If $t_4^4 + 1 = 0$, then two runs of the GBA yields an empty variety.

Now, let $t_2t_1 - t_4 = 0$. We do the same process. Start with $t_2 + 1 = 0$ and $t_2^2 + 1 = 0$. Both immediately give $p = 0$. If $t_2^4 + 1 = 0$, then two runs of the GBA returns an empty variety. Let $t_2^2 + t_2 + 1 = 0$. Two runs of the GBA yield $d_2 - 1 = 0 \implies d_4 - d_1 = 0$ and also that $t_2 = a_1$. This is enough to fill out the entire $S$-matrix. The missing elements at this point are $a_3$ and $a_4$. But Galois actions on the eigenvalue matrix fill that in for us. Note, that $s_{42} = a_3$ because $d_2 = 1$. Adding in the coefficients of the $s$-polynomial for $d_1$ returns $d_1^2 \pm d_1 - 1 = 0$. But since we can assume that $d_1 \geq 0$, we can assume $d_1$ is the golden ratio. Then, $\sigma$ we can tell sends $\phi$ to $-1/\phi$, where $\phi$ is the golden ratio. Observer that $\sigma(a_3) = \frac{-t_2}{\phi} \implies a_3 = \phi t_2$. And similarly $a_4 = \phi \bar{t}_2$. Without loss of generality,

$$S = \begin{pmatrix} 1 & \phi & 1 & 1 & \phi & \phi \\ \phi & -1 & \phi & \phi & -1 & -1 \\ 1 & \phi & t_2 & \bar{t}_2 & \phi t_2 & \phi \bar{t}_2 \\ 1 & \phi & \bar{t}_2 & t_2 & \phi \bar{t}_2 & \phi t_2 \\ \phi & -1 & \phi t_2 & \phi \bar{t}_2 & t_2 & \bar{t}_2 \\ \phi & -1 & \phi \bar{t}_2 & \phi t_2 & \bar{t}_2 & t_2 \end{pmatrix}$$

Since the submatrix associated with the pointed subcategory is invertible, it is a modular subcategory and hence is a product of the modular data of Fibonocci MTC and the $\mathbb{Z}_3$ MTC [9]. Descriptions of those categories can be found in [11] □

### 5.3. Self Dual.

**Lemma 28.** *Let $\mathcal{C}$ be a rank NI-SD-MTC. If $\mathcal{C}$ is a product category, then without loss of generality the Galois group is one of $\langle (01)(23) \rangle$, $\langle (012)(345) \rangle$, $\langle (01)(23)(45), (02)(13) \rangle$, or $\langle (012345) \rangle$.*

*Proof:* In [11], Rowell, Strong, and Wang classify all MTC's up to rank 4. If $\mathcal{C}$ is a product category, then it must be a product of a 2 dimenstional and a 3 dimensional category. It should be clear that both products must be SD, otherwise the $S$ matrix will have non-real entries and therefore $\mathcal{C}$ would be NSD. The 2 dimensional categories given in [11] are the Semion, and Fibonacci. Both are SD. The 3 dimensional SD categories are Ising, $(A_1, 2)$, and $(A_1, 5)_{\frac{1}{2}}$. The $S$ matrix is the same for both Ising and $(A_1, 2)$. This means the products will yield the same Galois groups. We only need to consider the products of Semion/Fibonacci and Ising/$(A_1, 5)_{\frac{1}{2}}$.

Semion $\otimes$ Ising gives $\langle (01)(23) \rangle$. Semion $\otimes (A_1, 5)_{\frac{1}{2}}$ gives $\langle (012)(345) \rangle$. Fibonacci $\otimes$ Ising gives $\langle (01)(23)(45), (02)(13) \rangle$. Lastly, Fibonacci $\otimes (A_1, 5)_{\frac{1}{2}}$ gives $\langle (012345) \rangle$. $\square$

**Theorem 29.** *There is no rank 6 NI-SD-MTC associated with $G = \langle (0123) \rangle$.*

*Proof:* Let $\sigma = (0123)$. Then,
$$S_{04} = \epsilon_1 \epsilon_4 S_{14} = \epsilon_1 \epsilon_2 S_{24} = \epsilon_1 \epsilon_2 \epsilon_3 \epsilon_4 S_{34} = \epsilon_1 \epsilon_2 \epsilon_3 \epsilon_0 S_{04} \implies \epsilon_1 \epsilon_2 \epsilon_3 \epsilon_0 = 1$$

But, $\prod \epsilon_i = -1 \implies \epsilon_4 \neq \epsilon_5 \implies S_{45} = \epsilon_4 \epsilon_5 S_{45} \implies S_{45} = 0$. WLOG, $\epsilon_4 = \epsilon_0$. This gives 4 sign choices,

1) $\epsilon_0 = \epsilon_1 = \epsilon_2 = \epsilon_3 = \epsilon_4 \neq \epsilon_5$  2) $\epsilon_0 = \epsilon_1 = \epsilon_4 \neq \epsilon_2 = \epsilon_3 = \epsilon_5$
3) $\epsilon_0 = \epsilon_2 = \epsilon_4 \neq \epsilon_1 = \epsilon_3 = \epsilon_5$  4) $\epsilon_0 = \epsilon_3 = \epsilon_4 \neq \epsilon_1 = \epsilon_2 = \epsilon_5$

$$S = \begin{pmatrix} 1 & d_1 & d_2 & d_3 & d_4 & d_5 \\ d_1 & \epsilon_1 \epsilon_2 d_2 & \epsilon_1 \epsilon_3 d_3 & \epsilon_0 \epsilon_1 & \epsilon_1 \epsilon_4 d_4 & \epsilon_1 \epsilon_5 d_5 \\ d_2 & \epsilon_1 \epsilon_3 d_3 & 1 & \epsilon_0 \epsilon_2 d_1 & \epsilon_1 \epsilon_2 d_4 & \epsilon_1 \epsilon_2 d_5 \\ d_3 & \epsilon_0 \epsilon_1 & \epsilon_0 \epsilon_2 d_1 & \epsilon_1 \epsilon_2 d_2 & \epsilon_0 \epsilon_4 d_4 & \epsilon_0 \epsilon_5 d_5 \\ d_4 & \epsilon_1 \epsilon_4 d_4 & \epsilon_1 \epsilon_2 d_4 & \epsilon_0 \epsilon_4 d_4 & a_1 & 0 \\ d_5 & \epsilon_1 \epsilon_5 d_5 & \epsilon_1 \epsilon_2 d_5 & \epsilon_0 \epsilon_5 d_5 & 0 & a_2 \end{pmatrix}$$

The FPcol cannot be column 4 or 5 since both have 0's. The other four columns share an orbit so we may assume that the FPcol is column 0 and $d_i > 0$ for all $i$. In choice 1, $0 = S_0 \cdot S_2$ is a sum of positive real numbers.

Note that Galois symmetry implies that if $t_3 - 1 = 0$ or $t_2 - t_3 = 0$, then $t_1, t_2, t_3 = 1$. Each sign choice requires the same deduction, namely after a series of runs of the GBA, the polynomial, $(t_3 - 1)(t_3 - t_2)(t_2 + 1)$ appears. The first two factors always occur at the same time as previously stated. They all lead to the same contradiction, therefore, $t_2 + 1 = 0$ and by Galois symmetry, $t_1 + t_3 = 0$ as well.

In all three sign choices, after the 4th run of the GBA, it's shown that $a_1 = a_2 = 0$. Then Orthogonality relations show $4d_4^2 = 4d_5^2 = D$. But, since $d_i > 0 \implies d_4 = d_5$. Therefore they should have the same Galois Conjugates in the same order in the eigenvalue matrix. In choice 2, the first Galois conjugates for $d_4$ and $d_5$ respectively are $\frac{d_4}{d_1}$ and $\frac{-d_5}{d_1}$. In choice 3, the corresponding ones are $\frac{-d_4}{d_1}$ and $\frac{d_5}{d_1}$. Finally in case 4, $\frac{-d_4}{d_1}$ and $\frac{d_5}{d_1}$. $\square$

**Theorem 30.** *There is no rank 6 NI-SD-MTC associated with $G = \langle (01234) \rangle$.*

*Proof:* If $t_i = t_j$ for any $i \neq j \implies t_i = 1$ for all $0 \leq i < 5$. If $t_5 = t_i$ for any $i \neq 5 \implies t_i = 1$ for all $i$. So, assume $t_i = 1$ for all $0 \leq i < 5$. Finding a contradiction for this is enough to show all $t_i$ are distinct. Let $\sigma = (01234)$. Also note that if $\sigma = (01234)$, $\prod \epsilon_i = 1$. Then,

$$S = \begin{pmatrix} 1 & d_1 & d_2 & d_3 & d_4 & d_5 \\ d_1 & \epsilon_1 \epsilon_2 d_2 & \epsilon_1 \epsilon_3 d_3 & \epsilon_1 \epsilon_4 d_4 & \epsilon_0 \epsilon_1 & \epsilon_1 \epsilon_5 d_5 \\ d_2 & \epsilon_1 \epsilon_3 d_3 & \epsilon_0 \epsilon_5 d_4 & \epsilon_3 \epsilon_5 & \epsilon_0 \epsilon_2 d_1 & \epsilon_1 \epsilon_2 d_5 \\ d_3 & \epsilon_1 \epsilon_4 d_4 & \epsilon_3 \epsilon_5 & \epsilon_1 \epsilon_5 d_1 & \epsilon_0 \epsilon_3 d_2 & \epsilon_0 \epsilon_4 d_5 \\ d_4 & \epsilon_0 \epsilon_1 & \epsilon_0 \epsilon_2 d_1 & \epsilon_0 \epsilon_3 d_2 & \epsilon_0 \epsilon_4 d_3 & \epsilon_0 \epsilon_5 d_5 \\ d_5 & \epsilon_1 \epsilon_5 d_5 & \epsilon_1 \epsilon_2 d_5 & \epsilon_0 \epsilon_4 d_5 & \epsilon_0 \epsilon_5 d_5 & z d_5 \end{pmatrix}$$

But after one run of the GBA, we get that $p = 0$. Therefore all $t_i$ are distinct. Since they are all distinct, this group cannot yield a modular category. By lemma 28 this is not a product category. So if there is any MTC, $\mathcal{C}$, with Galois group $\langle (01234) \rangle$, it will be nondegenerate and prime. According to [7] the only strong modular fusion algebra of dimension 6 described in table 12 is of level 9. But, $\text{Gal}(\mathbb{Q}_9/\mathbb{Q})$ doesn't have a subgroup of order 5. □

**Theorem 31.** *There is no rank 6 NI-SD-MTC associated with $G = \langle (01)(2345) \rangle$.*

*Proof:* Let $\sigma = (01)(2345)$. Then, $S_{02} = \epsilon_1\epsilon_2 S_{15} = \epsilon_0\epsilon_1\epsilon_2\epsilon_5 S_{04} = \epsilon_0\epsilon_2\epsilon_4\epsilon_5 S_{13} = \epsilon_2\epsilon_3\epsilon_4\epsilon_5 S_{02} \implies \epsilon_2\epsilon_3\epsilon_4\epsilon_5 = 1$. But, $\prod \epsilon_i = 1 \implies \epsilon_0 = \epsilon_1$. This gives,

$$S = \begin{pmatrix} 1 & d_1 & \cdots \\ d_1 & 1 & \\ d_2 & \epsilon_0\epsilon_3 d_3 & \\ d_3 & \epsilon_0\epsilon_3 d_2 & \\ \epsilon_2\epsilon_5 d_2 & \epsilon_0\epsilon_4 d_3 & \\ \epsilon_2\epsilon_3 d_3 & \epsilon_0\epsilon_2 d_2 & \cdots \end{pmatrix}$$

Now, orthogonality of the first two columns gives, $0 = 2d_1 + 4\epsilon_0\epsilon_3 d_2 d_3$, or $\frac{-\epsilon_0\epsilon_3}{2} = \frac{d_2 d_3}{d_1}$. But, $\frac{d_2 d_3}{d_1}$ is an integer. □

**Theorem 32.** *There is no rank 6 NI-SD-MTC associated with $G = \langle (01)(23), (23)(45) \rangle$.*

*Proof:* Let $\sigma = (01)(23)$ and $\tau = (23)(45)$. Then $\epsilon_0 = \epsilon_1$, $\epsilon_2 = \epsilon_3$, $\delta_2 = \delta_3$, and $\delta_4 = \delta_5$. But, $\prod \epsilon_i = \prod \delta_i = 1 \implies \epsilon_4 = \epsilon_5$ and $\delta_0 = \delta_1$.

Then,

$$S = \begin{pmatrix} 1 & d_1 & d_2 & \delta_0\delta_2 d_2 & d_4 & \delta_0\delta_4 d_4 \\ d_1 & 1 & \epsilon_0\epsilon_2\delta_0\delta_2 d_2 & \epsilon_0\epsilon_2 d_2 & \epsilon_0\epsilon_4 d_4 & \epsilon_0\epsilon_4\delta_0\delta_4 d_4 \\ d_2 & \epsilon_0\epsilon_2\delta_0\delta_2 d_2 & a_1 & a_2 & a_3 & \epsilon_2\epsilon_4\delta_2\delta_4 a_3 \\ \delta_0\delta_2 d_2 & \epsilon_0\epsilon_2 d_2 & a_2 & a_1 & \epsilon_2\epsilon_4 a_3 & \delta_2\delta_4 a_3 \\ d_4 & \epsilon_0\epsilon_4 d_4 & a_3 & \epsilon_2\epsilon_4 a_3 & a_4 & a_5 \\ \delta_0\delta_4 d_4 & \epsilon_0\epsilon_4\delta_0\delta_4 d_4 & \epsilon_2\epsilon_4\delta_2\delta_4 a_3 & \delta_2\delta_4 a_3 & a_5 & a_4 \end{pmatrix}$$

Orthoganaility of column 0 and column 1 give, $-1 = \frac{d_2^2}{d_1}\epsilon_0\epsilon_2\delta_0\delta_2 + \frac{d_4^2}{d_1}\epsilon_0\epsilon_4$. But $\frac{d_2^2}{d_1}$ and $\frac{d_4^2}{d_1}$ are integers with the same sign. Therefore, $\epsilon_0\epsilon_2\delta_0\delta_2 \neq \epsilon_0\epsilon_4 \implies \epsilon_2\epsilon_4 \neq \delta_0\delta_2$. This leaves 8 sign choices.

1) $\epsilon_0 = \epsilon_2 = \epsilon_4$   2) $\epsilon_0 = \epsilon_2 = \epsilon_4$   3) $\epsilon_0 = \epsilon_2 \neq \epsilon_4$   4) $\epsilon_0 = \epsilon_2 \neq \epsilon_4$
   $\delta_0 = \delta_4 \neq \delta_2$        $\delta_0 \neq \delta_2 = \delta_4$        $\delta_0 = \delta_2 = \delta_4$        $\delta_0 = \delta_2 \neq \delta_4$

5) $\epsilon_0 = \epsilon_4 \neq \epsilon_2$   6) $\epsilon_0 = \epsilon_4 \neq \epsilon_2$   7) $\epsilon_0 \neq \epsilon_2 = \epsilon_4$   8) $\epsilon_0 \neq \epsilon_2 = \epsilon_4$
   $\delta_0 = \delta_2 = \delta_4$        $\delta_0 = \delta_2 \neq \delta_4$        $\delta_0 = \delta_4 \neq \delta_2$        $\delta_0 \neq \delta_2 = \delta_4$

Choices 4, 8, 5 and 7 are relabelings of choices 1, 2, 3, and 6 respectively. Sign choice 6 does not have a column that can be all positive or all negative. Sign choices 1 and 2 have positive $d_2$ and negative $d_2$ in columns 0 through 3. Therefore none of those columns be entirely positive or entirely negative. The FPcol, must be either column 4 or 5. Either way, $a_3$ can't be 0. In sign choice 3 the reverse is true. Both positive and negative $a_3$ appear in columns 2 through 5. Therefore we can assume the FPcol is column 0.

In all cases we'll show that both 5 and 16 must divide $n$. We'll do this by showing that there are two elements of $\text{Gal}(\mathbb{Q}_n/\mathbb{Q})$ that act differently on the $t$-spectra via Galois symmetry. If only one of 5 or 16 divide $n$, then all elements of order 4 would act on the $t$-spectra the same via Galois symmetry. Recall that the action of $\sigma$ acts via $\sigma^2$.

Let $\mathcal{C}$ be a an MTC with $G = \langle (01)(23), (23)(45) \rangle$. If both 16 and 5 divide $n$ and give rise to $\sigma$'s that act differently, then both must have subrepresentations of dimension greater than one appearing in the decomposition and factorization of any modular representation of $\mathcal{C}$. The allowable representations of order 16 have dimension 3 or 6. If it's dimension 6, there isn't enough room in the $t$-spectra for both 5 and 16 to divide $n$. Therefore we assume $\rho = 2 \otimes 3$ where he 2-dimensional factor is order 5 and the 3 dimensional factor has order 16. This gives, $\{\zeta_5, \zeta_5^4\} \otimes \{\zeta_8, \zeta_{16}, -\zeta_{16}\} = \{\zeta_5\zeta_8, \zeta_5^4\zeta_8, \zeta_5\zeta_{16}, \zeta_5^2\zeta_{16}, -\zeta_5\zeta_{16}, -\zeta_5^4\zeta_{16}\}$. But the Galois symmetry of $\langle (01)(23), (23)(45)\rangle$ clearly doesn't lineup with this $t$-spectra.

To show that both 5 and 16 divide $n$, we must show that some pair liftings of $\sigma$, $\tau$, or $\sigma\tau$ act differently on the $t$-spectra via Galois symmetry. It will be enough then to show that there is at most 1 pair of identical $t_i$'s in $\{1, t_1, t_2, t_3, t_4, t_5\}$.

If $t_2 - 1 = 0 \implies t_3 - t_1 = 0$ by Galois symmetry of $\sigma$. But then, by Galois symmetry of $\tau$, $\gamma t_0$ is fixed by liftings of $\tau^2$. But since $t_0 = t_2$, $\gamma t_2$ must also be fixed by liftings of $\tau^2$. But liftings of $\tau^2$ are supposed to map $\gamma t_2$ to $\gamma t_3 \implies t_2 = t_3 \implies t_1 = 1$. Similar observations show that if $t_i = t_j$ where $i$ and $j$ are in different orbits then all $t_i$ in either orbit are identical. Therefore it's enough to show the following three cases lead to contradictions, $t_1 - 1 = t_2 - t3 = 0$, $t_1 - 1 = t_4 - t_5 = 0$, and $t_2 - t_3 = t_4 - t_5 = 0$.

The work for sign choices 1 and 2 are identical as $a_3 \neq 0$ in both choices. One run of the GBA assuming that $t_1 - 1 = 0$, yields $t_2 - t_3 = t_4 - t_5$. Two more runs of the GBA and it turns out that all $t_i = 1$ and therefore $N = 1$. This implies the Galois group is trivial. Next we assume that $t_2 - t_3 = t_4 - t_5 = 0$. One run yields that $t_1 - 1 = 0$. But that puts us back into the first case which we've already eliminated.

Sign choice 3 is a little more involved. If $t_1 - 1 = t_2 - t_3 = 0$, then one run of the GBA yields $a_3(t_4 - t_5)$. We first assume $t_4 - t_5 = 0$. After two runs of the GBA we can conclude that both $t_2$ and $t_4$ are 6th roots of unity (not necessarily primitive). Therefore $N$ divides 6. But $\text{Gal}(\mathbb{Q}_6/mathbbQ)$ does not have a proper subgroup of order 4.

At this point we know that $t_4 - t_5$ is not 0. This implies, either 5 or 16 but not both divide $n$. If 5 divides $n$, then we appeal to lemma 26. We need to find a possible $t$-spectra that has exactly pair of $\gamma t_i$ that are moved by $\tau^2$, but exactly two pair of identical $\gamma t_i$. Such a $t$-spectra doesn't exist. Therefore, 16 divides $n$ and $t_4 + t_5 = 0$. Adding that relation and $a_3 = 0$ into the GBA yields, $4d_4^2 - D = 4a_5^2 - D = a_4 + a_5 = 0$, i.e. $\pm d_4 = a_5 = -a_6$. But this immediately implies that columns 4 and 5 of the eigenvalue matrix are integer columns and therefore fixed by $G$.

Next we suppose $t_1 - 1 = t_4 - t_5 = 0$. Since we've already eliminated $t_1 - 1 = t_2 - t_3 = t_4 - t_5 = 0$ we can assume $t_2 - t_3 \neq 0$. It's still true that 5 doesn't divide $n$. So, 16 must and $t_2 + t_3 = 0$. Adding $t_1 - 1 = t_4 - t_5 = t_2 + t_3 = 0$ to the GBA and running it twice, yields a similar contradiction, but in columns 2 and 3.

Finally we enter the relations, $t_1 + 1 = t_2 - t_3 = t_4 - t_5 = 0$ and one run of the GBA yields $D = 0$. $\square$

**Theorem 33.** *There is no rank 6 NI-SD-MTC associated with $G = \langle (01)(23), (02)(13)\rangle$.*

*Proof:* Let $\sigma = (01)(23)$ and $\tau = (02)(13)$. Then $\epsilon_0 = \epsilon_1$, $\epsilon_2 = \epsilon_3$, $\delta_0 = \delta_2$, and $\delta_1 = \delta_3$. But, $\prod \epsilon_i = \prod \delta_i = 1 \implies \epsilon_4 = \epsilon_5$ and $\delta_4 = \delta_5$. Note that $S_{03} = \epsilon_0\epsilon_2 S_{12} = \delta_0\delta_1 S_{21} \implies \epsilon_0\epsilon_2 = \delta_0\delta_1$. Then,

$$S = \begin{pmatrix} 1 & d_1 & d_2 & d_3 & d_4 & d_5 \\ d_1 & 1 & \epsilon_0\epsilon_2 d_3 & \epsilon_0\epsilon_2 d_2 & \epsilon_0\epsilon_4 d_4 & \epsilon_0\epsilon_4 d_5 \\ d_2 & \epsilon_0\epsilon_2 d_3 & 1 & \epsilon_0\epsilon_2 d_1 & \delta_1\delta_4 d_4 & \delta_1\delta_4 d_5 \\ d_3 & \epsilon_0\epsilon_2 d_2 & \epsilon_0\epsilon_2 d_1 & 1 & \epsilon_0\epsilon_4\delta_0\delta_4 d_4 & \epsilon_0\epsilon_4\delta_0\delta_4 d_5 \\ d_4 & \epsilon_0\epsilon_4 d_4 & \delta_0\delta_4 d_4 & \epsilon_0\epsilon_4\delta_1\delta_4 d_4 & z_1 d_4 & z_2 d_5 \\ d_5 & \epsilon_0\epsilon_4 d_5 & \delta_0\delta_4 d_5 & \epsilon_0\epsilon_4\delta_1\delta_4 d_5 & z_3 d_4 & z_4 d_5 \end{pmatrix}$$

This leaves 8 cases.

1) $\epsilon_0 = \epsilon_2 = \epsilon_4$   2) $\epsilon_0 = \epsilon_2 = \epsilon_4$   3) $\epsilon_0 = \epsilon_2 \neq \epsilon_4$   4) $\epsilon_0 = \epsilon_2 \neq \epsilon_4$
   $\delta_0 = \delta_1 = \delta_4$        $\delta_0 = \delta_1 \neq \delta_4$     $\delta_0 = \delta_1 = \delta_4$      $\delta_0 = \delta_1 \neq \delta_4$
5) $\epsilon_0 = \epsilon_4 \neq \epsilon_2$   6) $\epsilon_0 = \epsilon_4 \neq \epsilon_2$   7) $\epsilon_0 \neq \epsilon_2 = \epsilon_4$   8) $\epsilon_0 \neq \epsilon_2 = \epsilon_4$
   $\delta_0 = \delta_4 \neq \delta_1$      $\delta_0 \neq \delta_1 = \delta_4$      $\delta_0 = \delta_4 \neq \delta_1$     $\delta_0 \neq \delta_1 = \delta_4$

Cases 1, 2, 3, and 4 have two columns of all positive entries. Cases 6 and 7 are relabelings of case 5. We can assume that the FPcol is column 0. If it were columns 4 or 5, then the FPdims would all be integers and therefore the category would be integral. Then, the remaining labels share an orbit. Thus we can assume $d_i \geq 1$ for all $i$.

For sign choice 5: One run of the GBA yields, $d_4^2 + d_5^2 + 2d_2 - 2d_3 = 0 \implies d_3 > d_2$. But it also yields, $(d_3+1)^2 - D = (d_2-1)((d_2-1)^2 - D) = 0$. If $d_2 - 1 = 0 \implies \frac{d_3}{d_1} = 1 \implies d_3 + d_1 = 0$. But, both are positive. So, $(d_3+1)^2 - D = (d_2-1)^2 - D \implies (d_2-1)^2 = (d_3+1)^2$. But since $d_i geq 1 \implies d_2 - 1 = d_3 + 1 \implies d_2 > d_3$, giving a contradiction.

For sign choice 8: One run of the GBA yields, $(d_2+1)^2 - D = (d_3+1)^2 - D = 0 \implies d_2 = d_3$. But $\frac{-d_3}{d_2} = -1$ is a Galois Conjugate of $d_1$, i.e. $d_1 = -1$. But $d_i \geq 1$, gives a contradiction.

**Theorem 34.** *There is no rank 6 NI-SD-MTC associated with $G = \langle (0123)(45) \rangle$.*

*Proof:* Let $\sigma = (0123)(45)$. Then consider, $S_{04} = \epsilon_1\epsilon_5 S_{15} = \epsilon_1\epsilon_5\epsilon_2\epsilon_4 S_{24} = \epsilon_1\epsilon_2\epsilon_4\epsilon_3 S_{35} = \epsilon_1\epsilon_2\epsilon_3\epsilon_0 S_{04} \implies \epsilon_0\epsilon_1\epsilon_2\epsilon_3 = 1$. But, $\prod \epsilon_i = 1 \implies \epsilon_4 = \epsilon_5$. If $\epsilon_0 = \epsilon_1 = \epsilon_2 = \epsilon_3$, then there will be at least two columns where every entry has the same sign. This leaves 6 sign choices.

1) $\epsilon_0 = \epsilon_1 = \epsilon_4 \neq \epsilon_2 = \epsilon_3$   2) $\epsilon_0 = \epsilon_2 = \epsilon_4 \neq \epsilon_1 = \epsilon_3$
3) $\epsilon_0 = \epsilon_3 = \epsilon_4 \neq \epsilon_1 = \epsilon_2$   4) $\epsilon_0 = \epsilon_1 \neq \epsilon_2 = \epsilon_3 = \epsilon_4$
5) $\epsilon_0 = \epsilon_2 \neq \epsilon_1 = \epsilon_3 = \epsilon_4$   6) $\epsilon_0 = \epsilon_3 \neq \epsilon_1 = \epsilon_2 = \epsilon_4$

$$S = \begin{pmatrix} 1 & d_1 & d_2 & d_3 & d_4 & d_5 \\ d_1 & \epsilon_1\epsilon_2 d_2 & \epsilon_1\epsilon_3 d_3 & \epsilon_0\epsilon_1 & \epsilon_1\epsilon_4 d_5 & \epsilon_1\epsilon_4 d_4 \\ d_2 & \epsilon_1\epsilon_3 d_3 & 1 & \epsilon_0\epsilon_2 d_1 & \epsilon_1\epsilon_2 d_4 & \epsilon_1\epsilon_2 d_5 \\ d_3 & \epsilon_0\epsilon_1 & \epsilon_0\epsilon_2 d_1 & \epsilon_1\epsilon_2 d_2 & \epsilon_0\epsilon_4 d_5 & \epsilon_0\epsilon_4 d_4 \\ d_4 & \epsilon_1\epsilon_4 d_5 & \epsilon_1\epsilon_2 d_4 & \epsilon_0\epsilon_4 d_5 & a_1 & a_2 \\ d_5 & \epsilon_1\epsilon_4 d_4 & \epsilon_1\epsilon_2 d_5 & \epsilon_0\epsilon_4 d_4 & a_2 & a_1 \end{pmatrix}$$

Then, sign choices 2 and 5 are relabelings of each other as are 3 and 6. We will show that all the $t_i$ are distinct. By Galois symmetry, if $t_1 = 1$ or $t_3 = 1$, then $t_1 = t_2 = t_3 = 1$. If $t_2 = 1$, then $t_1 = t_3$. If $t_4 = 1$ or $t_5 = 1$, then $t_1 = t_3$ and $t_2 = 1$. To show all $t_i$ distinct it suffices to show that $t_3, t_2 \neq 1$ and $t_4 \neq t_5$, due to Galois symmetry.

All four sign choices follow a similar pattern. A first run of the Gröbner basis algorithm yields two polynomials with three interesting factors, two of which are $t_3 - 1$ and $t_4 - t_5$. This is naturally helpful since we're already trying to eliminate both of those cases. Successfully doing so will give more relations when considering $t_2 - 1 = 0$. In all sign choices when assuming $t_4 - t_5 = 0$, $h_4 = 4d_5^2 + a_1^2 + a_2^2 - D$. This relation implies that $d_4^2 = d_5^2$. We can assume the FP column is column 0 and therefore that $d_i \geq 1 \implies d_4 = d_5$.

In sign choices 1, 2, and 4 under the same assumption, the 7th run of the GBA, yields $t_5^2 - t2$ and $t_1 t_3 - t_2$ as relations. Recall from Galois symmetry, $\sigma^2(\gamma t_i) = \gamma t_{\sigma(i)}$. Since, $t_0 = 1$, $\sigma^2(\gamma) = t_1$ and $\sigma^4(\gamma) = t_2$. Note, $\sigma^4(\gamma t_2) = t_0 \implies \sigma^4(t_2) = \frac{1}{t_2}$. But, $\sigma^4(t_2) = \sigma^8(t_5) = t_5 \implies \frac{1}{t_2} = t_5 \implies t_2^3 = 1$. But if this is true, then $N = 1$ or 6. In either case, $\mathbb{Z}_N^\times$ does not have a cyclic subgroup of order 4. In all sign choices all other contradictions necessary to eliminate a subcase involved showing a non-zero element (or a product) such as $p$ or $d_4 d_5$ was 0. Most of the time this was directly given as a result of the GBA. But the final contradiction in each sign choice came from the GBA yielding $a_2 = 0$. Then, by the Galois action on the eigenvalue matrix, $a_1 = 0$. Then $0 = S_4 \cdot S_5 = 4 d_4 d_5$.

But neither $d_4$ nor $d_5$ can be zero. Therefore all $t_i$'s are distinct. This is enough due to [7] since it is also not a product category nor listed on his table $\square$

**Theorem 35.** *There is no rank 6 NI-SD-MTC associated with $G = \langle (012), (345) \rangle$.*

*Proof:* Let $\sigma = (012)$. We know there is only one sign choice, $\epsilon_0 = \epsilon_2 \neq \epsilon_1 = \epsilon_3$ because of theorem 12. Let $\tau = (345)$. Consider, $S_{01} = \delta_0 \delta_1 S_{01} \implies \delta_0 = \delta_1$. Similarly, $S_{02} = \delta_0 \delta_2 \implies \delta_0 = \delta_1 = \delta_2$. Then, $\prod \delta_i = 1 \implies \delta_0 \delta_3 \delta_4 \delta_5 = 1$. This gives,

$$S = \begin{pmatrix} 1 & d_1 & d_2 & d_3 & \delta_0\delta_3 d_3 & \delta_3\delta_4 d_3 \\ d_1 & -d_2 & -1 & d_3 & \delta_0\delta_3 d_3 & \delta_3\delta_4 d_3 \\ d_2 & -1 & d_1 & -d_3 & -\delta_0\delta_3 d_3 & -\delta_3\delta_4 d_3 \\ d_3 & d_3 & -d_3 & a_1 & a_2 & a_3 \\ \delta_0\delta_3 d_3 & \delta_0\delta_3 d_3 & -\delta_0\delta_3 d_3 & a_2 & \delta_4\delta_5 a_3 & \delta_3\delta_4 a_1 \\ \delta_3\delta_4 d_3 & \delta_3\delta_4 d_3 & -\delta_3\delta_4 d_3 & a_3 & \delta_3\delta_4 a_1 & \delta_3\delta_5 a_2 \end{pmatrix}$$

Instead of going into sign choices of the $\delta$ function, we include the restrictions in the GBA, $\delta_i^2 = 1$ and $\delta_0 \delta_3 \delta_4 \delta_5 = 1$. Then one run of the GBA yields, $t_3 D = 0$ and neither can be. This is a contradictions. $\square$

**Theorem 36.** *There is no rank 6 NI-SD-MTC associated with $G = \langle (01)(23)(45), (24)(35) \rangle$.*

*Proof:* Let $\sigma = (01)(23)(45)$ and $\tau = (24)(35)$. Then, $S_{02} = \epsilon_1 \epsilon_2 S_{13}$, $S_{20} = \epsilon_3 \epsilon_0 S_{31}$, $S_{04} = \epsilon_1 \epsilon_4 S_{15}$, $S_{40} = \epsilon_5 \epsilon_0 S_{51}$, and $\prod \epsilon_i = -1 \implies \epsilon_0 \epsilon_1 = \epsilon_2 \epsilon_3 = \epsilon_4 \epsilon_5 = -1$. Also note that $\delta_2 = \delta_4$, $\delta_3 = \delta_5$, and $\prod \delta_i = 1 \implies \delta_0 = \delta_1$.

$$S_{02} = \epsilon_1 \epsilon_2 S_{13} = \epsilon_1 \epsilon_2 \delta_0 \delta_3 S_{15}$$
$$S_{02} = \delta_0 \delta_2 S_{04} = \delta_0 \delta_2 \epsilon_1 \epsilon_4 S_{15} \implies \delta_2 \delta_3 = \epsilon_2 \epsilon_4$$
$$S_{03} = \epsilon_1 \epsilon_3 S_{12} = \epsilon_1 \epsilon_3 \delta_0 \delta_2 S_{14}$$
$$S_{03} = \delta_0 \delta_3 S_{05} = \delta_0 \delta_3 \epsilon_1 \epsilon_5 S_{14} \implies \delta_2 \delta_3 = \epsilon_3 \epsilon_5$$

This gives 8 sign choices:
1) $\epsilon_0 = \epsilon_2 = \epsilon_4 \neq \epsilon_1 = \epsilon_3 = \epsilon_5$
   $\delta_0 = \delta_2 = \delta_4$
2) $\epsilon_0 = \epsilon_2 = \epsilon_4 \neq \epsilon_1 = \epsilon_3 = \epsilon_5$
   $\delta_0 \neq \delta_2 = \delta_3$
3) $\epsilon_0 = \epsilon_2 = \epsilon_5 \neq \epsilon_1 = \epsilon_3 = \epsilon_4$
   $\delta_0 = \delta_2 \neq \delta_3$
4) $\epsilon_0 = \epsilon_2 = \epsilon_5 \neq \epsilon_1 = \epsilon_3 = \epsilon_4$
   $\delta_0 = \delta_3 \neq \delta_2$
5) $\epsilon_0 = \epsilon_3 = \epsilon_4 \neq \epsilon_1 = \epsilon_2 = \epsilon_5$
   $\delta_0 = \delta_2 \neq \delta_3$
6) $\epsilon_0 = \epsilon_3 = \epsilon_4 \neq \epsilon_1 = \epsilon_2 = \epsilon_5$
   $\delta_0 = \delta_3 \neq \delta_2$
7) $\epsilon_0 = \epsilon_3 = \epsilon_5 \neq \epsilon_1 = \epsilon_2 = \epsilon_4$
   $\delta_0 = \delta_2 = \delta_3$
8) $\epsilon_0 = \epsilon_3 = \epsilon_5 \neq \epsilon_1 = \epsilon_2 = \epsilon_4$
   $\delta_0 \neq \delta_2 = \delta_3$

$$S = \begin{pmatrix} 1 & d_1 & d_2 & d_3 & \delta_0\delta_2 d_2 & \delta_0\delta_3 d_3 \\ d_1 & -1 & \epsilon_1\epsilon_3 d_3 & \epsilon_1\epsilon_2 d_2 & \epsilon_1\epsilon_5\delta_0\delta_3 d_3 & \epsilon_1\epsilon_4\delta_0\delta_2 d_2 \\ d_2 & \epsilon_1\epsilon_3 d_3 & a_1 & a_2 & a_3 & a_4 \\ d_3 & \epsilon_1\epsilon_2 d_2 & a_2 & -a_1 & \epsilon_3\epsilon_5 a_4 & \epsilon_3\epsilon_4 a_3 \\ \delta_0\delta_2 d_2 & \epsilon_1\epsilon_5\delta_0\delta_3 d_3 & a_3 & \epsilon_3\epsilon_5 a_4 & a_1 & \epsilon_3\epsilon_5 a_2 \\ \delta_0\delta_3 d_3 & \epsilon_1\epsilon_4\delta_0\delta_2 d_2 & a_4 & \epsilon_3\epsilon_4 a_3 & \epsilon_3\epsilon_5 a_2 & -a_1 \end{pmatrix}$$

For relabelings, switch labels 2 and 3 and also switch labels 4 and 5. Doing so gives, sign choices 7, 8, 6, and 5 as relabelings of 1, 2, 3, and 4 respectively. We eliminate this group by showing that for all sign choices the $t_i$'s are distinct. Without loss of generality, we only need to show, $t_1, t_2, t_3 \neq 1$ and $t_3, t_4, t_5 \neq t_2$. Unless specified, all sign choices yielded the same results from the GBA.

Assume, $t_1 = 1$. After one run of the GBA, we get that $t_2 - t_4 = t_3 - t_5 = t_4 + t_5 + 1 = 0$. This means that $t_4$ and $t_5$ are complex conjugates and the real part is -1/2. The only complex units with real part -1/2 are the two primitive third roots of unity. Therefore, $N = 3$. The Galois group of Ga,$(\mathbb{Q}(\zeta_3)/\mathbb{Q}) \cong \mathbb{Z}/2\mathbb{Z}$, which doesn't have a subgroup of order 4.

If, $t_2 = 1$. Then by Galois symmetry, $t_1 - t_3 = t_2 - t_4 = t_3 - t_5 = 0$. Similarly if $t_3 = 1$, then, $t_1 - t_2 = t_2 - t_4 = t_3 - t_5 = 0$. After one run of the GBA, we get $t_5^2 - 1 = 0$ and $t_4^2 - 1$ respectively. Both then imply that $N = 1$ or 2. But, $\mathbb{Q}(\zeta_i) = \mathbb{Q}$ for $i = 1, 2$. Thus the Galois group is trivial and doesn't have a subgroup of order 4.

If $t_2 = t_3$ then by Galois symmetry, $t_4 = t_5$. After one run of the GBA, we get that $pD = 0$.

Note that by Galois symmetry, $t_2 = t_4 \iff t_3 = t_4$ and $t_2 = t_5 \iff t_3 = t_4$. If $t_2 = t_5$, then after one run of the GBA, we get that $t_1 - 1 = 0$, which we've already shown leads to a contradiction.

The final case, $t_2 = t_4$, is the most involved. We can assume that $1, t_1, t_2,$ and $t_3$ are all distinct. Then, using Galois symmetry, we can conclude that all elements of $\text{Gal}(\mathbb{Q}_n/\mathbb{Q})$ have order 1, 2, or 4 and any element of order 4 is a lifting of $\sigma$ and all liftings of $\sigma$ have order 4. This means that $n$ divides 240, and that 5 or 16 divide $n$. Both cannot divide $n$, otherwise the elements of order 4 projections onto $\text{Gal}(\mathbb{Q}_n/\mathbb{Q})$ would have different images in $\text{Sym}_6$.

Suppose 5 divides $n$. We must now appeal to the modular representations. We know that 5 divides $n$ so we begin there. Under the assumption we have exactly two pairs of identical elements in the $t$-spectra. Both of these pairs must be moved by all the liftings of $\sigma^2$. The final two elements in the $t$-spectra are also moved to each other under the same liftings (i.e. every element in the $t$-spectra will be moved).

These two pair also mean that the modular representation can only be decomposed into a direct sum of at most 3 subrepresentations. This is due to the restriction of overlapping $t$-spectra. All the $\mathbb{Z}/5\mathbb{Z}$ representations of dimension 5 or 6 give exactly one pair of identical elements in the $t$-spectra. Any combination including a $\mathbb{Z}/5\mathbb{Z}$ representation of dimension 3 will have two elements fixed.

From lemma 26 there are two possible $t$-spectra, $\{\alpha_1\} \otimes \{\zeta_5, \zeta_5^4, \zeta_5, \zeta_5^4, \zeta_5^2, \zeta_5^3\}$ or $\{\alpha_1\zeta_5, \alpha_1\zeta_5^4, \alpha_1\zeta_5, \alpha_1\zeta_5^4, \alpha_2, \zeta_5, \alpha_2\zeta_5^4\}$, $\zeta_5$ is any primitive 5th root of unity and $\alpha_1, \alpha_2$ are some 24th roots of unity. Recall that the $t$-spectra is just $\gamma \otimes T$-spectra for some root of unity $\gamma$. Since the $T$-spectra contains 1, then the $t$-spectra contains $\gamma$. Then it's clear that no matter which element of the first possible $t$-spectra is $\gamma$, all factors of $\alpha_1$ will be removed, leaving $T^5 = 1$ or $N = 5$. However, $\text{Gal}(\mathbb{Q}_5/\mathbb{Q})$ is a cyclic group of order 4 and therefore does not contain a subgroup isomorphic to $\mathbb{Z}/2\mathbb{Z} \oplus \mathbb{Z}/2\mathbb{Z}$.

Consider, $\{\alpha_1\zeta_5, \alpha_1\zeta_5^4, \alpha_1\zeta_5, \alpha_1\zeta_5^4, \alpha_2, \zeta_5, \alpha_2\zeta_5^4\}$. Under the current assumptions, the 1 in the $T$-spectra is unique, i.e. without loss of generality $\gamma = \alpha_2\zeta_5$ and the $T$-spectra $= \{1, \zeta_5^3, \alpha, \alpha\zeta_5^3, \alpha, \alpha\zeta_5^3\}$, where $\alpha = \alpha_1\alpha_2$. This gives us the following new relations for the GBA, $t_1^4 + t_1^3 + t_1^2 + t_1 + 1 = 0$ and $(t_3 - t_1t_2)(t_2 - t_1t_3) = 0$.

In all 4 sign choices of the $S$-matrix, this breaks down into two cases, $t_3 - t_1t_2 = 0$ and $t_2 - t_1t_3 = 0$. Both cases require 2 runs of the GBA each. The first run yields $a_1 + a_3 + 1 = Dt_5^2 + p(1 + t_4 + t - t_5 - 1) = 0$ (as the case or the sign choice varies, the signs may change or $t_4$ and $t_5$ may switch places). The second run yields, $t_4^2 + 2a_3 - t_4 + 1 = 0$ (as the case or the sign choice changes, the sign of $2a_3$ may change or $t_4$ may be replaced with $t_5$).

Consider, $t_4^2 - t_4 = -2a_3 - 1$. This implies that $t_4^2 - t_4$ is real because the entries of $S$ are real in all self-dual modular categories. Thus, $t_4^2 = t_4$ or $t_4^2 = -\bar{t}_4 \implies t_4 = 1$ or $-t_4\bar{t}_4 = t_4^3 = -1$. If $t_4 = 1$, then $a_3 = -1/2$ (or positive 1/2 depending on case and sign choice). But then $a_3$ would not be an algebraic integer. If $t_4^3 = -1 \implies t_4^2 - t_4 + 1 = 0 \implies a_3 = 0$. But $a_3 = 0 \implies a_1, a_2, a_4 = 0$ as well. Then, orthogonality relations imply that $d_2^2 + d_3^2 - D = 0$ and $1 + d_1^2 + 2d_2^2 + 2d_3^2 - D = 0$. Together, these imply that $1 + d_1^2 + d_2^2 + d_3^2 = 0$. This is a sum of positive real numbers. It cannot be equal to 0.

Now suppose 16 divides $n \implies n = 16$ or 48. From the lemma proved earlier, this implies $t_2 + t_3 = 0$. This provides an extra relation for the GBA. After two runs, each sign choice yielded $Dd_2d_3 = 0$.

Therefore in any sign choice $t_i \neq t_j$ for $i \neq j$. $\square$

**Theorem 37.** *There is no rank 6 NI-SD-MTC associated with $\langle(01),(2345)\rangle$ or $\langle(01),(23)(45),(24)(35)\rangle$.*

*Proof:* Let $\sigma = (01)$. From theorem 9, for any $2 \leq i, j \leq 5$ such that $\epsilon_i \neq \epsilon_j$, $S_{ij} = 0$, and such a pair, $i, j$ exist. Due to relabeling we can assume the following two cases, $\epsilon_5 \neq \epsilon_i$ for $2 \leq i < 5$ and $\epsilon_0 = \epsilon_1 = \epsilon_2 = \epsilon_3 \neq \epsilon_4 = \epsilon_5$. In the first case, let $\epsilon = \epsilon_0\epsilon_2$. This gives the following two matrices,

$$S1 = \begin{pmatrix} 1 & d_1 & d_2 & d_3 & d_4 & d_5 \\ d_1 & 1 & \epsilon d_2 & \epsilon d_3 & \epsilon d_4 & -\epsilon d_5 \\ d_2 & \epsilon d_2 & a_1 & a_2 & a_3 & 0 \\ d_3 & \epsilon d_3 & a_2 & a_4 & a_5 & 0 \\ d_4 & \epsilon d_4 & a_3 & a_5 & a_6 & 0 \\ d_5 & -\epsilon d_5 & 0 & 0 & 0 & a_7 \end{pmatrix}$$

$$S2 = \begin{pmatrix} 1 & d_1 & d_2 & d_3 & d_4 & d_5 \\ d_1 & 1 & d_2 & d_3 & -d_4 & -d_5 \\ d_2 & d_2 & a_1 & a_2 & 0 & 0 \\ d_3 & d_3 & a_2 & a_3 & 0 & 0 \\ d_4 & -d_4 & 0 & 0 & a_4 & a_5 \\ d_5 & -d_5 & 0 & 0 & a_5 & a_6 \end{pmatrix}$$

A priori, there are 5 possible Galois groups that contain (01). We'll show that $S1$ and $S2$ only allow for 1 group each aside from $G = \langle(01)\rangle$. Note in both matrices each row starting with row 2 has a 0. Consider the groups, $\langle(01),(2345)\rangle$ and $\langle(01),(23)(45),(24)(35)\rangle$. Both of those show the orbit of column 2 is the columns 2 through 5. This would force $a_i = 0$ for all $i$ in both $S1$ and $S2$. In particular $S_2 \cdot S_3 = 2d_2d_3 = 0$ by the orthogonality relations ($S_2$ and $S_3$ are the second and third column of either $S1$ or $S2$). But neither $d_2$ nor $d_3$ are 0. So, neither $\langle(01),(2345)\rangle$ nor $\langle(01),(23)(45),(24)(35)\rangle$ is possible. Indeed, in the case of $S1$ the only possible groups are $\langle(01)\rangle$ and $\langle(01),(234)\rangle$. And in the case of $S2$ the only possible groups are $\langle(01)\rangle$ and $\langle(01),(23)(45)\rangle$. $\square$

**Theorem 38.** *There is no rank 6 NI-SD-MTC associated with $G = \langle(01)(23)(45)\rangle$.*

*Proof:* This particular group is a subgroup of two Galois groups that do indeed have associated categories. So this group is much more difficult to prove no category is associated with it. Let $\sigma = (01)(23)(45)$. Then, $S_{02} = \epsilon_0\epsilon_1\epsilon_2\epsilon_3 S_{02}$ and $S_{04} = \epsilon_0\epsilon_1\epsilon_4\epsilon_5 S_{04}$. But then, $\epsilon_0\epsilon_1\epsilon_2\epsilon_3 = \epsilon_0\epsilon_1\epsilon_4\epsilon_5 = 1 \implies \epsilon_0\epsilon_1 = \epsilon_2\epsilon_3 = \epsilon_4\epsilon_5$. Note that, $\prod \epsilon_i = -1 \implies \epsilon_0 \neq \epsilon_1$, $\epsilon_2 \neq \epsilon_3$, and $\epsilon_4 \neq \epsilon_5$. This gives 4 sign choices.

1) $\epsilon_0 = \epsilon_2 = \epsilon_4 \neq \epsilon_1 = \epsilon_3 = \epsilon_5$   2) $\epsilon_0 = \epsilon_2 = \epsilon_5 \neq \epsilon_1 = \epsilon_3 = \epsilon_4$
3) $\epsilon_0 = \epsilon_3 = \epsilon_4 \neq \epsilon_1 = \epsilon_2 = \epsilon_5$   4) $\epsilon_0 = \epsilon_3 = \epsilon_5 \neq \epsilon_1 = \epsilon_2 = \epsilon_4$

But all cases are relabelings of each other. We approach this group with representation theory from the beginning. At the same time we focus on what separates this group from the larger groups that do have associated categories, the $s$-polynomials of the various entries. Combined with the variety of cases from representation theory, this is enough to eliminate this group.

First we establish some preliminary restrictions on the $T$-spectra. We want to say there at most two pair of identical $t_i$'s. To do so we rule 4 and 6 identical $t_i$ and three pair of identical $t_i$. If $t_i = 1$ for all $i$, one run of the GBA yields $p = 0$. If $t_1 = t_2 = t_3 = 1$, $t_1 = t_4 = t_5 = 1$, or $t_2 = t_3 = t_4 = t_5$, then the GBA yields $pD = 0$. If $t_1 = 1$, $t_2 = t_3$, and $t_4 = t_5$, then the GBA yields $p = 0$. If $t_1 = 1$, $t_2 = t_4$, and $t_3 = t_5$ or $t_1 = 1$, $t_2 = t_5$, and $t_3 = t_4$, then the GBA yields $t_4 + t_5 + 1 = 0$. But that implies that $t_4$ and $t_5$ are primitive $3^{rd}$ roots of unity and $N = 3$. Then, $\text{Gal}(\mathbb{Q}_N/\mathbb{Q})$ does not have a proper subgroup of order 2.

If $t_3 = 1$, $t_2 = t_1$, and $t_4 = t_5$ then the GBA yields, $b_2(t_2 + t_5 + 1) = 0$. If $t_2 + t_5 + 1 = 0$ then $N = 3$. Therefore $b_2 = 0 \implies b_1 = 0$. The GBA also yields, $d_3b_2t_2 - 2d_5c_2t_5 + d_3b_2 = 0$. Since $b_2 = 0$, this simplifies to $2d_5c_2t_5 = 0 \implies c_2 = 0 \implies c_1 = 0$. But then, $D = S_4 \cdot S_4 = d_4^2 + d_5^2 \implies 1 + d_1^2 + d_2^2 + d_3^2 = 0$, a sum of positive real numbers. Similar contradiction is found if $t_2 = 1$, $t_3 = t_1$, and $t_4 = t_5$.

If $t_5 = 1$, $t_4 = t_1$, and $t_2 = t_3$, then the GBA yields $a_2(t_3 + t_4 + 1) = 0$. Again, $t_3 + t_4 + 1 \implies N = 3 \implies a_2 = 0 \implies a_1 = 0$. The GBA also yields $a_2c_1t_4 - a_1c_2t_4 + c_2pt_4 - a_2c_1 + a_1c_2 + c_2p = 0 \implies c_2p(t_4 + 1) = 0$. But, $c_2 = 0$ yields a similar contradiction to the previous one, therefore, $t_4 + 1 = 0$. The GBA also yielded, $b_1t_4 + b_1 + t_3 = 0 \implies t_3 = 0$. A similar contradiction is found if $t_4 = 1$, $t_5 = t_1$, and $t_2 = t_3$. Now we can assume that there are at most two pair of identical $t_i$'s. But also Galois symmetry now implies that either 5 or 16 divide $n$ but not both.

Suppose 16 divides $n$. We can then assume that if $t_1 \neq 1$, then $t_1 = -1$, if $t_2 \neq t_3$, then $t_2 = -t_3$, and if $t_4 \neq t_5$ then $t_4 = -t_5$. We begin by considering all the possible ways there could be 2 pair of identical $t_i$ and then all the ways there could be 1 pair of identical $t_i$. These will all lead to contradictions. That will leave us with 16 divides $n$ implies that all the $t_i$ are distinct. This is enough to say that the any modular representation cannot be decomposed as a direct sum of subrepresentations 5.3. But, [7] says that the smallest representation of order 16 has dimension 6, i.e not only does $\rho$ not decompose, it doesn't factor either (1 dimensional factors are still allowed). Eholzer calls this simple and nondegenerate. In his table 12 in [7] there is only one simple nondegenerate representation that is admissible of dimension 6, and the order of its $t$ matrix is 9.

We used relations to define the coefficients of the $s$-polynomials. These were coded in and labeled as $f$'s as other $s$-polynomial coefficients have been in other groups. But we also used relations that were labeled $g$'s. These defined the elements of the eigenvalue matrix $s$. A first GBA was run to give relations labeled with $k$'s. These are immediate consequences of the orthogonality relations when put in terms of both the eigenvalue matrix and the $S$-matrix entries.

If $t_1 - 1 = t_2 - t_3 = t_4 + t_5 = 0$ then, after 3 runs of the GBA and using the $g$'s and $k$'s described above, we determined that $d_3$ was an integer and therefore, $d_1d_3 + d_2 = 0$. One more run and we were able to conclude that $d_1 = d_2$, and eventually that $d_1$ is an int, which implies $d_1 = \frac{-1}{d_1} \implies d_1^2 = -1$. But $S$ is a real matrix. A similar contradiction was found in similar steps when $t_1 - 1 = t_2 + t_3 = t_4 - t_5 = 0$.

If $t_2 - 1 = t_3 - t_1 = t_1 + 1 = t_4 + t_5 = 0$, then $(a_1 + 1)(p^2 - D) = b_2(p^2 - D) = 0$. If $c_2$ and $p^2 - D$ aren't 0, then, $b_2 = a_1 + 1 = 0 \implies b_1 = 0$. This eventually leads again to $d_1^2 = -1$ and therefore $p^- D = 0$. This shows $p$ is real. Eventually the GBA yields $2d_2d_3t_5 - d_1pt_5 + 2d_4d_5 - 2d_1t_5 = 0 \implies t_5$ is real, because $p$ is also real. Similar steps lead to similar contradictions for $t_3 - 1 = t_2 - t_1 = t_1 + 1 = t_4 + t_5 = 0$, $t_4 - 1 = t_5 - t_1 = t_1 + 1 = t_2 + t_3 = 0$, and $t_5 - 1 = t_4 - t_1 = t_1 + 1 = t_2 + t_3 = 0$.

If $t_2 - t_3 = t_4 - t_5 = t_1 + 1 = 0$, the GBA yields an empty variety.

If $t_2 - t_4 = t_3 - t_5 = t_1 + 1$, then the GBA yields 12 relations of the form $(p^2 - D)(\text{XXXX}) = 0$. If $p^2 - D = 0$, then after one run of the GBA we deduce that $D - 4 = 0$. But this yields, $(t_5 + 1)(d_3^2 + d_5^2) = 0$. If $t_5 + 1 = 0$, then $N = 2$ and $d_3^2 + d_5^2$ is a sum of positive real numbers and therefore can't be 0. So, the 12 XXXX relations must be 0. The GBA immediately yields $d_5 D = 0$. Similar steps lead to a similar contradiction if $t_2 - t_5 = t_3 - t_4 = t_1 + 1 = 0$.

If $t_1 - 1 = t_2 + t_3 = t_4 + t_5 = 0$, then the GBA yields $p = 0$. If $t_1 + 1 = t_2 - t_3 = t_4 + t_5 = 0$, then the GBA yields $pD = 0$. If $t_1 + 1 = t_2 + t_3 = t_4 - t_5 = 0$, then the GBA yields $Dd_5 = 0$. Therefore if 16 divides $n$, all the $t_i$ are distinct.

Suppose 5 divides $n$. We will work down the list from lemma 26. In each possible $t$-spectra, we'll have subcases depending on what element of the $t$-spectra is $\gamma$. Most of the time without loss of generality, there will be two choices for $\gamma$. The choice of $\gamma$ will determine what $t_1$ is. The Galois symmetry will determine what $t_2$ and $t_3$ are, but not which is which, i.e. we will have even more subcases. Luckily, all the subcases will be handled in similar manners. Also as a reminder, we can relabel the indexes by switching 2 and 4 and switching 3 and 5. This will help keep the number of subcases somewhat manageable.

At a certain point we will use [10]. We will show that $N = 5$ and thus, 2 is relatively prime to $N$. In point (v) of definition 2, we see that $\nu_2(i) = \pm 1$ if $i^* = i$. Our supposed category is self-dual, so this holds. But NG's theorem actually forces $\nu_2(i) = 1$, since you can't have negative coefficients in a direct sum. The orbit of 0 is $\{0, 1\}$. We know that $X_0^2 = X_0 \neq \sum_i X_i$. Therefore, we can conclude that $X_1^2 = \sum_i X_i$. Passing through to the $\dim(X_1^2)$, implies that $d_1^2 = 1 + d_1 + d_2 + d_3 + d_4 + d_5$, under the correct assumptions. There will be other times that we use (vii) from definition 2 as well. This group represents the first time either of these results are needed in the paper. The $f$'s as before define the coefficients of the $s$-polynomials of particular elements of $S$ or $s$.

If the $t$-spectra is $\alpha \otimes \{1, 1, \zeta_5, \zeta_5^4, \zeta_5^2, \zeta_5^3\} \implies$ the $T$-spectra is $\{1, 1, \zeta_5, \zeta_5^4, \zeta_5^2, \zeta_5^3\}$ or $\{1, \zeta_5^3, \zeta_5^4, \zeta_5^4, \zeta_5, \zeta_5^2\}$.

Suppose the $T$-spectra is $\{1, 1, \zeta_5, \zeta_5^4, \zeta_5^2, \zeta_5^3\}$. Then, $t1 - 1 = t_2^4 - t_3 = t_4^4 - t_5 = t_2^4 + t_3^3 + t_2^2 + t_2 + 1 = 0$. But also, $(t_2^2 - t_4)(t_2^2 - t_5) = 0$. Now suppose $t_2^2 - t_4 = 0$. In the second run of the GBA we add all 16 $f$'s and it yields, $z_{11}^2 = 2z_{15}^2$. But then clearly $z_{11}$ is not an integer. As similar contradiction is found if $t_2^2 - t_5 = 0$.

Suppose the $T$-spectra is $\{1, \zeta_5^3, \zeta_5^4, \zeta_5^4, \zeta_5, \zeta_5^2\}$. Then, $t_2 - t_3 = t_1^4 + t_1^3 + t_1^2 + t_1 + 1 = 0$. We also get, $(t_1 t_4 - 1)(t_1 t_5 - 1) = (t_2 t_5 - 1)(t_2 t_4 - 1) = (t_5^2 - t_4)(t_4^2 - t_5) = 0$. In all three products, the first terms happen at the same time and the second terms happen at the same time.

Suppose, $t_1 t_4 - 1 = t_2 t_5 - 1 = t_5^2 - t_4 = 0$. Then, after a few runs where we deduce that $b_1 - t_3 - t_5 - 1 = 0$ and that $2t_3 + t_4 + 2$ is not 0 (both $t_3$ and $t_4$ are $5^{th}$ roots), we show that $z_6$ which is the algebraic norm of $d_4, d_5$ is 2. But (vii) of definition 2 says that the prime ideals of $d_i$ are also prime ideals of $D$ and $N$. But the ideal generated by 2 is clearly not a divisor of the ideal generated by 5. A similar contradiction is found if $t_1 t_5 - 1 = t_2 t_4 - 1 = t_4^2 - t_5 = 0$.

Thus (2) of lemma 26 is not possible. Suppose the $t$-spectra is $\{\zeta_5, \zeta_5^4, \zeta_5, \zeta_5^4, \zeta_5^2, \zeta_5^3\}$. Then the $T$-spectra is $\{1, \zeta_5^3, 1, \zeta_5^3, \zeta_5, \zeta_5^2\}$ or $\{1, \zeta_5^3, \zeta_5, \zeta_5^2, \zeta_5, \zeta_5^2\}$.

Suppose the $T$-spectra is $\{\zeta_5, \zeta_5^4, \zeta_5, \zeta_5^4, \zeta_5^2, \zeta_5^3\}$. Then, $t_1^4 + t_1^3 + t_1^2 + t_1 + 1 = 0$. Here we have two cases to deal with, $(t_2 - 1)(t_3 - 1) = 0$ and $(t_4^2 - t_5)(t_5^2 - t_4) = 0$. But these two cases are independent. Suppose $t_2 - 1 = 0 \implies t_3 - t_1 = 0$. And also suppose $t_4^2 - t_5 = 0$.

The first run of the GBA provides enough new relations that the the second run would produce an ideal of 1 dimension. In some sense this means that there is 1 degree of freedom and that if we had one more useful relation we could then tell the GBA to attempt to solve for a specific variable. That's exactly what we do. Adding NG's relation, is enough to get a 0 dimension ideal. So, attempting to solve for one variable will work. Until this point we've been using the default monomial order of Macaulay2. The default ordering is a graded ordering by degree and then within a degree is ordered by a lexigraphical ordering based on the variables chosen for the polynomial ring at the beginning. To make use of NG's relation, we'll define a new polynomial ring and use a new monomial order called Eliminate. We'll tell Macaulay to try to eliminate the first 17 variables (there are 18 in our polynomial ring). Doing so with a 0 dimensional ideal guarantees the first polynomial will be a polynomial of the $18^{th}$ variable, in this case $c_2$. The GBA does indeed yield a polynomial in $c_2$, it has a handful of factors, most of which have leading coefficient that is not 1. Since $c_2$ is an algebraic integer, none of those irreducible factors can be 0. This leaves, $c_2(c_2-1) = 0$. If $c_2 = 0 \implies c_1 = 0$. Together with another relation from the first run of the GBA also force $b_1 = b_2 = 0$. But then $D = d_4^2 + d_5^2 \implies 0 = 1 + d_1^2 + d_2^2 + d_3^2$. If $c_2 - 1 = 0$, then $d_5^2 - 2 = 0$. But 2 doesn't divide 5. Similar contradictions are found for the other subcases.

Suppose the $T$-spectra is $\{1, \zeta_5^3, \zeta_5, \zeta_5^2, \zeta_5, \zeta_5^2\}$. Then $t_1^4 + t_1^3 + t_1^2 + t_1 + 1 = t_2 t_3 - t_1 = 0$. We again have two independent cases, $(t_2 - t_4)(t_2 - t_5) = 0$ and $(t_2 t_1 - 1)(t_3 t_1 - 1) = 0$. Suppose $t_2 - t_4 = 0 \implies t_3 - t_5 = 0$. Also suppose that $t_3 t_1 - 1 = 0 \implies t_2^2 - t_3 = 0$. After two runs of the GBA (and adding the $f$'s into the second run), we get that $z_3 + z_6 + 3 = 0$. But $z_3$ and $z_6$ are the algebraic norms of $d_2$ and $d_4$ respectively. They must share their prime divisors with $N$, i.e. they must be powers of 5 (0 is an allowed power). But no two powers of 5 will ever differ by 3. A similar contradiction is found in the other subcases.

Suppose the $t$-spectra is $\alpha \otimes \{\alpha_1 \zeta_5, \alpha_1 \zeta_5^4, \alpha_1 \zeta_5, \alpha_1 \zeta_5^4, \alpha_2 \zeta_5, \alpha_2 \zeta_5^4\}$. Then without loss of generality, the $T$-spectra is $\{1, \zeta_5, 1, \zeta_5, \alpha, \alpha\zeta_5\}$ or
$\{1, \zeta_5, \alpha, \alpha\zeta_5, \alpha, \alpha\zeta_5\}$ where $\alpha$ is some primitive $2^{nd}$, $3^{rd}$, $4^{th}$, or $8^{th}$ root of unity.

Suppose the $T$-spectra is $\{1, \zeta_5, 1, \zeta_5, \alpha, \alpha\zeta_5\}$. Then $t_1^4 + t_1^3 + t_1^2 + t_1 + 1 = 0$. Once again, we have two independent cases, $(t_2 - 1)(t_3 - 1) = 0$ and $(t_4 t_1 - t_5)(t_5 t_1 - t_4) = 0$. Even after those 4 cases, each will have 4 subcases based on $\alpha$. In all cases when $\alpha$ is a $3^{rd}$ or $4^{th}$ root of unity we get a straightforward contradiction. The contradiction associated to $\alpha$ being a primitive $8^{th}$ root is more interesting. After 3 runs of the GBA we get that $\alpha$ is real or $t_5^3 - 1 = 0$ (possibly $t_4^3 - 1 = 0$ depending on subcases). But $\alpha$ is a primitive $8^{th}$ root and can't be real and $t_5$ is a primitive $40^{th}$ root.

Finally if $\alpha = -1$ we have some deviation in the subcases. Three of them have $\sigma(d_1) = -1/d_1$ and $+1/d_1$. But one subcase, $t_2 - 1 = t_4 t_1 - t_5 = \alpha + 1 = 0$, doesn't immediately have that contradiction. It actually allows for the modular data to be factored as a tensor product.

The modular data factors a tensor product of Fibonacci category, see a description in [11] and an integral matrix. Since Fibonacci is known to be modular the other matrix is an integer matrix, but there are no self dual integral MTC of rank three.

Suppose the $T$-spectra is $\{1, \zeta_5, \alpha, \alpha\zeta_5, \alpha, \alpha\zeta_5\}$. Then $t_1^4 + t_1^3 + t_1^2 + t_1 + 1 = 0$. One more time we have 2 independent cases, $(t_2 - t_4)(t_2 - t_5) = 0$ and $(t_2 t_1 - t_3)(t_3 t_1 - t_2) = 0$. Each subcase requires 3 runs of the GBA and yields something like, $t_4^2 - a_1 - b_1 - t_4 - 1 = 0$. But then $t_4^2 - t_4$ is real. This means that $t_4 - 1 = 0$ or $t_4$ is a primitive $6^{th}$ root. But $t_4$ is $\alpha$ (in this particular subcase, in others we might switch $t_5$ and $t_4$). And $\alpha$ cannot be a primitive $6^{th}$ root of unity.

Suppose the $t$-spectra is $\alpha \otimes \{1, 1, \alpha_1, \alpha_2, \zeta_5, \zeta_5^4\}$. But this doesn't fit the Galois symmetry of $\langle (01)(23)(45) \rangle$. Take the $\alpha_1$ and $\alpha_2$ pair. They differ, so according to the Galois symmetry, if $\Sigma$ is a lifting $\sigma = (01)(23)(45)$, then $\Sigma^2(\alpha_1) = \alpha_2$. But since $\alpha_1$ is some $24^{th}$ root of unity, $\Sigma^2(\alpha_1) = \alpha_1$ for all $\Sigma$. $\square$

## 6. Conclusion and Future Work

In this paper, we showed the efficacy of using a computational approach based on the admissibility criteria in [3] to classify low rank modular categories by giving a partial classification of rank 6 modular categories. Using the Galois structure and action on the eigenvalue matrix, we constructed all the possible abelian subgroups of $\text{Sym}_6$ for two pair of non-self-dual simple objects, one pair of non-self-dual objects, and the case where all simple objects are self-dual.

Our approach begins by fixing a subgroup of $\text{Sym}_r$ and then building all the possible modular data based on the Galois symmetry of the $S$ matrix, $S_{ij} = \epsilon(\sigma(i))\epsilon(j) S_{\sigma(i)\sigma^{-1}(j)}$. Here $\sigma$ is an element of the image of $\text{Gal}(\mathbb{Q}(S)/\mathbb{Q})$ in $\text{Sym}_r$ and $\epsilon$ is a sign function dependent on the choice of $\sigma$. Then, we run an initial Gröbner basis calculation, (GBA), using the orthogonality and twist relations found in definition 2.

The output of a (GBA) is a basis for the ideal of polynomials generated by the input. Given ideal conditions a basis would lead to solving the modular data. It isn't a surprise that we do not begin with ideal conditions. Much of the facts in the admissibility criteria are not algebraic, i.e. they cannot be written in terms of a polynomial. For example, $S_{0j} > 0$ for all $j$, but there is no polynomial that has only the positive real numbers as a solution. Therefore we often ran multiple GBA calculations. In-between each, we factored the output and looked at the results that factored. If a polynomial in the output did indeed factor, we used the admissibility criteria to try to determine if one factor was never allowed to be zero. The most common examples were in the form of $S_{0j} h_i$ for some $j$. Then, $h_i$ was added to the next GBA calculation. The process stops when a contradiction is found (i.e. no category is associated with the given pair of $S$ and $T$), or the modular data is solved. The most common contradiction came in from some element of the output forcing one of the following elements to be 0; $d_j, p,$ or $D$.

We showed that a rank 6, non-integral, and non-self-dual category is isomorphic to a product of two modular categories. We did this by first eliminating the following groups as described above, for two pair of non-self-dual objects: $\langle (01), (23), (45) \rangle$, $\langle (01)(23), (23)(45) \rangle$, $\langle (01), (23)(45) \rangle$, $\langle (01)(2435) \rangle$, $\langle (01), (2435) \rangle$, and $\langle (01), (23)(45), (24)(35) \rangle$, for one pair of non-self-dual objects: $\langle (01), (23), (45) \rangle$, $\langle (012), (45) \rangle$, $\langle (01), (45) \rangle$, $\langle (01)(23), (02)(13), (45) \rangle$, and $\langle (0213), (45) \rangle$. The remaining two groups, $\langle (01)(23), (45) \rangle$ (one pair of NSD objects) and $\langle (01)(24)(35), (23)(45) \rangle$ (two pair of NSD objects) required the use of representation theory. Representation theory allowed us to eliminate $\langle (01)(23), (45) \rangle$ and prove that all modular data from the group $\langle (01)(24)(35), (23)(45) \rangle$ is a tensor product of two modular categories of rank 2 and 3.

Similarly we eliminated the following groups assuming all simple objects are self-dual, $\langle (0123) \rangle$, $\langle (01234) \rangle$, $\langle (01)(23), (02)(13) \rangle$, $\langle (0123)(45) \rangle$, $\langle (012), (345) \rangle$, $\langle (01), (2345) \rangle$, $\langle (01)(2345) \rangle$, $\langle (01), (23)(45), (24)(35) \rangle$, $\langle (01)(23)(45), (24)(35) \rangle$, $\langle (01)(23)(45) \rangle$, or $\langle (01)(23), (23)(45) \rangle$. It is known that following groups do have categories associated with them, $\langle (012) \rangle$, $\langle (01)(23) \rangle$, $\langle (012)(345) \rangle$, $\langle (01)(23)(45), (02)(13) \rangle$, and $\langle (012345) \rangle$. It is unknown but we conjecture that the following groups do not have a modular category associated to them, $\langle (01) \rangle$, $\langle (01), (234) \rangle$, and $\langle (01), (23)(45) \rangle$.

The approach described in this paper has a natural extension to classifying modular categories of other low rank. Specifically I intend to complete the classification of rank 6 in the self-dual case and then attempt a classification of non-integral rank 7 modular categories. It is also possible that this approach can be adapted to classify low rank fusion categories with slightly different structure such as super modular categories.

## References

[1] B. Bakalov and A. K. Jr. *Lectures on Tensor Categories and Modular Functors*. Vol. 21. University Lecture Series. Amer. Math. Soc., 2001.

[2] J. de Boer and J. Goeree. "Markov Traces and $\Pi_1$ Factors in Conformal Field Theory". *Comm. Math. Phys.* 139.2 (1991), pp. 267–304.

[3] P. Bruillard, S.-H. Ng, E. Rowell, and Z. Wang. "On Classification of Modular Categories by Rank". *Internat. Math. Res. Not.* 2016 (2016), pp. 7546–7588.

[4] P. Bruillard, S.-H. Ng, E. Rowell, and Z. Wang. "Rank-finiteness for modular categories". *J. Amer. Math. Soc.* 29.3 (2016), pp. 857–881.

[5] C. Dong, X. Lin, and S.-H. Ng. "Congruence Property in Conformal Field Theory". *Algebra Number Theory* 9.9 (2015), pp. 2121–2166.

[6] W. Eholzer. "Fusion algebras induced by representations of the modular group". *Internat. J. Modern Phys.* 8 (1993), pp. 3495–3507.

[7]  W. Eholzer. "On the Classification of Modular Fusion Algebras". *Comm. Math. Physics.* 172.3 (1995), pp. 623–659.
[8]  P. Etingof, S. Gelaki, D. Nikshych, and V. Ostrik. *Tensor Categories*. Vol. 205. Mathematical Surveys and Monograhs. Amer. Math. Soc., 2015.
[9]  M. Müger. "On Fusion Categories". *Proc. London Math. Soc.* 87.2 (2003), pp. 291–308.
[10] S.-H. Ng. Conference Presentation. 2017 West Coast Lie Theory Workshop. 6-03-2017.
[11] E. Rowell, R. Stong, and Z. Wang. "On Classification of Modular Tensor Categories". *Comm. Math. Phys.* 292.2 (2009), pp. 343–389.
[12] V. Turaev. *Quantum Invariants of Knots and 3-Manifolds*. Vol. De Gruyter Studies in Mathematics. Walter de Gruyter, 1994.